\documentclass[graybox]{svmult}
\usepackage{type1cm}        
\usepackage{makeidx}         
\usepackage{graphicx}        
\usepackage{multicol}        
\usepackage[bottom]{footmisc}

\usepackage{newtxtext}       %
\usepackage[varvw]{newtxmath}       

\makeindex             

\begin{document}

\title*{Constructive semigroups with apartness - a state of the art}
\titlerunning{Constructive semigroup with apartness}
\author{Melanija Mitrovi\'c, Mahouton Norbert Hounkonnou and Paula Catarino}
\institute{Melanija Mitrovi\'c \at Faculty of Mechanical Engineering, CAM-FMEN, University of Ni\v s, Serbia, \email{melanija.mitrovic@masfak.ni.ac.rs}
\and Mahouton Norbert Hounkonnou \at University of Abomey-Calavi, 072 B.P. 50 Cotonou, Republic of Benin, \email{norbert.hounkonnou@cipma.uac.bj}
\and Paula Catarino \at Department of Mathematics, University of Tr\'{a}s-os-Montes e Alto Douro, Vila Real, Portugal \email{pcatarin@utad.pt}}

\maketitle
{
This chapter aims  to provide a clear and understandable picture of
constructive semigroups with apartness in Bishop's style of
constructive mathematics, \textbf{BISH}. Our theory is partly inspired 
by the classical case,  but it is distinguished from it in two significant aspects:
we use \emph{intuitionistic logic}  rather than classical throughout;
our work is based on the notion of \emph{apartness }
(between elements of the set, and, later, between elements and its
subsets). Following Heyting, at least initially, classical semigroup theory is
seen as a guide that helps us to develop the constructive theory of
semigroups with apartness. To have a structure, we need a set, a relation, and rules establishing
how we will put them together. Working within classical or intuitionistic logic, in order to
 analyze algebraic structures, it is necessary to start with study on
 sets and ordered sets, relational systems, etc. A
comparative analysis between presented classical and constructive
results is also a part of this chapter. All proofs can be found
in the Appendix.}


\section{Introduction}\label{ch-msa-1}

The development of powerful computer systems draws attention to the
intuitive notion of an effective procedure,  and of computation in
general. This, in turn, stimulates the development of constructive
algebra and its possible applications. \emph{Apartness}, the second
most important fundamental notion developed in constructive
mathematics, shows up in computer science. Inspired by results
obtained in interactive theorem proving the approach of formal
verifications (more in Mitrovi\'c, Hounkonnou and Baroni, 2021),  a new constructive algebraic theory known as
the \textbf{theory of semigroups with apartness} was developed by
 Mitrovi\'c and co-authors: Hounkonnou, Baroni, Crvenkovi\'c,
 Romano, Silvestrov - see (Crvenkovi\'c  {\it et al}, 2013), (Crvenkovi\'c  {\it et al}, 2016),
(Mitrovi\'c  {\it et al}, 2021),  (Mitrovi\'c and Silvestrov, 2020), (Mitrovi\'c   {\it et al}, 2019).
 We did it, we developed
the theory of constructive semigroups with apartness, and we are
still working on its development. But, is it Art?

A descriptive definition of a semigroup with apartness includes two
main parts:
\begin{itemize}
\item  the notion of a certain classical algebraic structure is
straightforwardly adopted;
\item  a structure is equipped with an apartness with a standard operation
respecting that apartness.
\end{itemize}
This chapter aims  to provide a clear and understandable picture of
constructive semigroups with apartness in Bishop's style of
constructive mathematics, \textbf{BISH}, both to (classical)
algebraists and those who apply algebraic knowledge.

Our theory is partly inspired by the classical case,  but it is
distinguished from it in two significant aspects:
\begin{enumerate}
\item we use \emph{\textbf{intuitionistic logic }} rather than
classical throught;
\item our work is based on the notion of \textbf{apartness }
(between elements of the set, and, later, between elements and its
subsets).
\end{enumerate}

``\emph{Intuitionistic logic} can be understood as the logic of
scientific research (rather  positivistically conceived) while on
the other hand the classical logic is the logic of ontological
thought,'' Grzegorczyk (1964). Intuitionistic logic
serves as a foundation to constructive mathematics. This means that
mathematics is done without use of the so-called The Law  of the
excluded middle,\textbf{LEM},  and without the use of
nonconstructive existence proofs.

It was A. Heyting (1956) who wrote in his
famous book \emph{Intuitionism, an Introduction} , that
``Intuitionism can only flourish, if mathematicians, working in
different fields, become actively interested in it and make
contributions to it. [...] In order to build up a definite branch of
intuitionistic mathematics, it is necessary in the first place to
have a thorough knowledge of the corresponding branch of classical
mathematics, and in the second place to know by experience where the
intuitionistic pitfalls lie.'' So, we present two points of view on
constructive semigroups with apartness:
\begin{itemize}
\item the classical (\textbf{CLASS}), which plays a useful role as intuition
guides and to at least link with the presentations given in
constructive one;
\item the constructive (\textbf{BISH}), comprising a binary system with apartness
which  satisfies a number of extra conditions: well known axioms of
apartness, and the operations have to be strongly extensional.
\end{itemize}

In
{\color{blue}https://webpages.uidaho.edu/rwells/techdocs/Biological}
it can be found that ``The
Bourbaki set out to discover the 'roots' of mathematics - something
that was true of mathematics in general. They found it - or so they
tell us - subsisting in three basic "mother structures" upon which
all of mathematics depends. \emph{These structures are not reducible
one to another. This just means no one of them can be derived from
the other two. But by making specifications within one structure, or
by combining two or more of these structures, everything else in
mathematics can be generated.} These three basic structures are
called \textbf{algebraic structure}, \textbf{order structure} and
topological structure.''

In what follows some results concerning special ordered and
algebraic structures - the \emph{binary} ones, e.g., the ordered
structures with a single binary relation, the algebraic structures
with a single binary operation and the ordered algebraic structures
with a single binary operation and a single binary relation -- such
as, certain types of ordered sets, semigroups and ordered semigroups
-- within classical and constructive (Bishop's) settings will be
presented.

 ``To have a structure, we need a set, a relation, and rules establishing
how we will put them together,''
{\color{blue}https://webpages.uidaho.edu/rwells/techdocs/Biological}.
That is, working
within classical or intuitionistic logic, in order to
 analyze algebraic structures it is necessary to start with study on
 sets and ordered sets, relational systems, etc. Selected results on
 classical ones are given in Subsection~\ref{ch-msa-21}, while
 constructive ones are given in Subsection~\ref{ch-msa-31}.

Following Heyting, at least initially, classical semigroup theory is
seen as a guide that helps us to develop the constructive theory of
semigroups with apartness. The classical background is given in
Section~\ref{ch-msa-2}. All in all, material presented in Section
~\ref{ch-msa-2} is  broad rather than deep, and it is not intended
to be comprehensive.  It is heavily based  the treatments in other
standard on set and semigroup theory textbooks. The main novelty is
in the selection and arrangement of material. Development of an
appropriate constructive order theory for sets and semigroups with
apartness has been one of the main objectives of our work. Most of
the results we have published to the date are given in
Section~\ref{ch-msa-3}. Applications and possible applications are
given in Subsection~\ref{ch-msa-23} for the classical, and in
Subsection~\ref{ch-msa-33} for the constructive case.
Subsection~\ref{ch-msa-24} and Subsection~\ref{ch-msa-34} contain
very short historical overviews and bibliographic notes, which give
sources and suggestions for further reading for both cases. A
comparative analysis between presented classical and constructive
results are part of Section~\ref{ch-msa-4}. All proofs can be found
in the Appendix.

Results presented in Section~\ref{ch-msa-3} are based on the ones
published in (Crvenkovi\'c  {\it et al}, 2013), (Crvenkovi\'c  {\it et al}, 2016),
(Mitrovi\'c  { \it et al}, 2021),  (Mitrovi\'c and Silvestrov, 2020), (Mitrovi\'c  {\it et al}, 2019),
(Darp\"{o} and Mitrovi\'c, arXiv:2103.07105), (Romano, 2005),  (Romano, 2007),  and
shortly, on the work in progress (Mitrovi\'c and Hounkonnou).

\section{The \textbf{CLASS} case}\label{ch-msa-2}

Can the whole of MODERN ALGEBRA be described in a couple of
sentences? Paraphrasing J. Wo$\acute{z}$ny (2018)  yes it
can; it has been designed to be elegantly simple: the story starts
with \textbf{sets} (collections of objects) and \textbf{relations}
on them  and proceeds to the concept of a \textbf{semigroup}; each
new concept is based on the previous ones, and, ultimately, the
whole multistory edifice \emph{rests on the sparse foundation of
sets and relations}. \textbf{Semigroups} serve as the building
blocks for the structures comprising the subject which is today
called \textbf{modern algebra}. In fact, as it is written by Uday S.
Reddy
({\color{blue}https://www.cs.bham.ac.uk/ udr/notes/semigroups.pdf})
``Semigroups are everywhere. Groups are
semigroups with a unit and inverses. Rings are ``double
semigroups:'' an inner semigroup that is required to be a
commutative group and an outer semigroup that does not have any
additional requirements. But the outer semigroup must distribute
over the inner one. We can weaken the requirement that the inner
semigroup be a group, i.e., no need for inverses, and we have
semirings. Semilattices are semigroups whose multiplication is
idempotent. Lattices are ``double semigroups'' again and we may or
may not require distributivity.''

\subsection{Sets and relations}\label{ch-msa-21}

It is well-known that all mathematical theories deal with sets in
one way or another, i.e  nearly every mathematical object of
interest is a set of some kind. In the words of Halmos (1960),
``General set theory is pretty trivial stuff
really, but, if you want to be a mathematician, you need some, and
here it is; read it, absorb it, and forget it.''

Many of the fundamental concepts of mathematics can be described in
terms of relations. The notion of an order plays an important role
throughout mathematics. A pure order theory is concerned with a
single undefined binary relation $\rho$. This relation is assumed to
have certain properties (such as, for example, reflexivity,
transitivity, symmetry, antisymmetry), the most basic of which leads
to the concept of \emph{quasiorder}, a reflexive and transitive
relation. A quasiorder plays a central role throughout this short
exposition. Combining quasiorder with the various properties of
symmetry  yields  the notions of an equivalence, a symmetric
quasiorder, and a (partial) order, an antisymmetric quasiorder.

Equivalences  permeate all of mathematics. An equivalence relation
allows one to connect those elements of a set that have a particular
property in common. In universal algebra, the formulation of mapping
images is one of the principal tools used to manipulate sets. In the
study of mapping images of a set a lot of help comes from the notion
of a quotient set, which captures all mapping images, at least up to
bijection. On the other hand, mapping is the concept which goes hand
in hand with equivalences. Thus concepts of equivalence, quotient
set and mappings are closely related. Knowing that the equivalence
$\varepsilon$ on a set $S$  is the kernel of the quotient map from
$S$ onto $S/\varepsilon$, we can treat equivalence relations on $S$
as kernels of mappings  with $S$ as the domain. The relationship
between quotients, mappings and equivalences is described by the
celebrated isomorphism theorems, which are a general and important
foundational part of universal algebra.

On the other hand, as Davey and Priestley (1990, 2002) wrote
  ``order, order, order - it permeates
mathematics, and everyday life'', as well ``to such an extent that
we take it as granted.''

If a set carries an order relation as well as an equivalence
relation, the important question arises: does the given order induce
an order on the equivalence classes? It is a well-known fact that
this will not be possible in general. Of course, there are
conditions guaranteeing that the quotient set inherits the ordering.

Applications of mathematics in other sciences often take the form of
a network of relations between certain objects, which justifies H.
Poincar\`{e} ``Mathematicians do not study objects, but relations
between objects,'' (quoted in Newman, 1956).

\subsubsection{Basic concepts and important examples}\label{ch-msa-211}

Mappings are written on the left, that is if we have a function $f$
that takes an input $x$, its output is written $f(x)$. If two maps
are composed, they are written right-to-left. The composition of two
mappings $f$ and $g$  is written $f\circ g$ or simply $fg$, and we
have $fg(x)= f(g(x))$. In connection with denoting  subsets, the
symbol $\subseteq$ is used to denote containment with possible
equality, while $\subset$ is used to denote strict containment.

Let $S$ and $T$ be sets.  A \emph{mapping} $f$ from $S$ to $T$,
denoted by $f: S \rightarrow T$, is a subset of $S\times T$ such
that for any element $x\in S$ there is precisely one element $y\in
Y$ for which $(x,y)\in f$, i.e.
$$(\forall  x,y\in S) \, x= y  \Rightarrow  f(x)= f(y).$$
Instead of $(x,y)\in f$, we usually write $y=f(x)$. Two mappings
$f,g:S \rightarrow T$ are \emph{equal} if they are equal as subsets
of $S\times T$, that is $f\, = \,  g \ \Leftrightarrow \ (\forall
x\in S )\, f(x) = g(x)$.

For the set of all mappings from $S$ to $T$ the notation $M(S,T)$ or
$T^S$ is used, and $\mathcal{T}_S$  for $M(S,S)$. The elements of
$\mathcal{T}_S$ are called \emph{transformations} of $S$.

A mapping $f$ is:

$\cdot$ \emph{surjective} or \emph{onto}:\quad $(\forall y\in
T)\,(\exists x\in S)\ (y= f(x))$;

$\cdot$ \emph{injective} or \emph{one-one}:
\begin{eqnarray*}
f(x)= f(y)\ &\Rightarrow& \ x= y\\
&\text{or},&\ \text{equivalently},\\
 x\neq y\  &\Rightarrow& \ f(x)\neq f(y);
\end{eqnarray*}

$\cdot$ \emph{bijective}: surjective and injective.

The subset of all bijective transformations in $\mathcal{T}_S$ is
denoted by $\mathcal{G}_S$, its elements are called
\emph{permutations} or \emph{bijections} of $S$. In particular, the
\emph{identity function} $i\iota_S : S \rightarrow S$, $i\iota_S(x)
= x$, is always a bijection, the \emph{trivial bijection}, of S. The
term 'function' is fundamental to almost all areas of mathematics,
but it has become traditional to use terms 'mapping' and
'transformation' in algebra.

Let $S$ be a set. A subset $\rho$ of $S\times S$, or, equivalently,
a property applicable to elements of $S\times S$, is called a
\emph{binary relation on} $S$. The empty subset $\emptyset$ of
$S\times S$ is included among the binary relations on $S$. Other
special binary relations worth to be mentioned are the universal
relation $S\times S$, and the \emph{equality} relation or
\emph{diagonal} $\triangle_S = \{(x,x)\,:\, x\in S\}$. Binary
relations on a set $S$ are often called \emph{homogenous}. The set
of all binary relation on $S$ is denoted by $\pmb{\mathcal{B}_S}$.

\subsubsection{Quotient sets. Homomorphism and isomorphism
theorems}\label{ch-msa-212}


In general, there are many properties that a binary relation may
satisfy on a given set. For example, the relation
$\rho\in\mathcal{B}_S$  might be:

\medskip

\begin{itemize}
\item[(R)] \quad  reflexive: $(x,x)\in\rho$ \quad ($\Leftrightarrow
\triangle_S\subseteq \rho$)
\item[(S)] \quad symmetric : $(x,y)\in\rho \ \Rightarrow\
(y,x)\in\rho$ \quad ($\Leftrightarrow \rho^{-1}\subseteq \rho
\Leftrightarrow \rho^{-1} = \rho)$
\item[(AS)] \quad  antisymmetric: $(x,y)\in\rho \wedge  (y,x)\in\rho  \ \Rightarrow\
x=y$ \quad  ($\Leftrightarrow \rho\cap \rho^{-1}\subseteq
\triangle_S$)
\item[(T)] \quad  transitive: $(x,y)\in\rho \wedge  (y,z)\in\rho  \ \Rightarrow\ (x,z)\in
\rho$
\end{itemize}

The most basic concept leads to the
\emph{\textbf{\textbf{quasiorders}}}, reflexive and transitive
relations with the fundamental concepts being introduced whenever
possible in their natural properties.
\begin{itemize}
\item[$\bullet$] \quad
 equivalence -  symmetric quasiorder;
\item[$\bullet$] \quad (partial) order - antisymmetric quasiorder.
\end{itemize}

The concept of an \emph{equivalence} is an extremely important one
and plays a central role in mathematics.  Classifying objects
according to some property is a frequent procedure in many fields.
Grouping elements in ``a company'' so that elements  in each group
are of the same prescribed property as performed by equivalence
relations, and the classification gives the corresponding quotient
sets. Thus, abstract algebra can show us how to identify objects
with the same properties properly - we have to switch to a quotient
structure (a technique applicable, for example, to abstract data
type theory).

Let us remember that we can define a surjective mapping starting
from equivalences, i.e. we have the following lemma.
\begin{lemma}\label{msart-11o}
Let $\varepsilon$ be an equivalence on $S$.  The mapping  $\pi : S
\rightarrow S/\varepsilon$ defined by $\pi (x)= x\varepsilon$, $x\in
S$, is a surjective mapping.
\end{lemma}

The mapping $\pi$ from the previous lemma is called the
\emph{natural mapping associated with an equivalence relation}
$\varepsilon$ and is sometimes written as $\pi_{\varepsilon}$ or as
$\pi_{nat}$. As usual $x\varepsilon$ denotes the $\varepsilon$-class
of an element $x\in S$.

On the other hand, any mapping $f$ between sets $S$ and $T$,
$f:S\rightarrow T$, gives rise to the equivalence on $S$
$$ker\,f = \{(x,y)\in S\times S\,:\, f(x)=f(y)\},$$
 called the \emph{kernel of mapping}  \ $f$ or \emph{the equivalence
relation induced by} $f$.

\emph{The First isomorphism theorem} for sets follows.

\begin{theorem}\label{msart-11}
Let $f:S \rightarrow T$ be a mapping between two sets $S$ and $T$.
Then the following statements are true:
\begin{itemize}
\item[\emph{(i)}] \quad  The mapping $\varphi : S/ker\,f \rightarrow T$
defined by $\varphi(x(ker\,f)) = f(x)$ is a unique injective mapping
such that $f=\varphi\circ\pi_{nat}$.
\item[\emph{(ii)}] \quad If $f$ is surjective, then $\varphi$ is a
bijection, i.e. we have $T \cong S/ker\,f$.
\item[\emph{(iii)}] \quad $\varphi$ is surjective if and only if $f$ is
surjective.
\end{itemize}
\end{theorem}

The next theorem concerns  a more general situation, and will be
frequently used.
\begin{theorem}\label{msart-12}
Let $\varepsilon$ be an equivalence on $S$. Let $f:S \rightarrow T$
be a mapping between two sets $S$ and $T$ such that
$\varepsilon\subseteq ker\,f$. Then there is a unique mapping
$\varphi: S/\varepsilon \rightarrow T$ defined by
$\varphi(x\varepsilon) = f(x)$ such that
$f=\varphi\circ\pi_{\varepsilon}$.
\end{theorem}

\subsubsection{Ordered sets}\label{ch-msa-213}

In many disciplines, sets which are equipped with a structure are
investigated.  The ordered pair $(S,\rho)$ where $S$ is a given set
and $\rho$ is a particular (binary)  relation defined on it is
called a (binary) \emph{relational structure}. A relational system
$(S,\rho)$ is called \emph{quasiordered} or \emph{ordered}  if
$\rho$ is  a quasiorder or an order  on $S$.  In what follows, we
will automatically assume that a property of a (quasi)order is
defined in the same way as the property of the (quasi)ordered set,
and vice versa.

 Relational structures play
an important role both in mathematics and in applications since
formal description of a real system often involve relations. For
these considerations we often ask about a certain factorization of a
relational system $(S,\rho)$ because of its importance in  enabling
us to introduce the method of abstraction on $S$. More precisely, if
$\varepsilon$ is an equivalence relation on $S$, we often ask about
a factor relation $\Theta_\rho=\rho/\varepsilon$ on the factor set
$S/\varepsilon$ such that the relational factor set
$(S/\varepsilon,\rho/\varepsilon)$ shares some ``good'' properties
of $(S,\rho)$.

Let $(S,\rho)$ and $(T,\sigma)$ be two relational systems. A mapping
$f:S\rightarrow T$ is called
\begin{itemize}
\item[--] \emph{isotone}  if  $(x,y)\in \rho \ \Rightarrow \ (f (x),
f(y))\in \sigma $;
\item[--]  \emph{reverse isotone} if $(f (x), f(y))\in \sigma
\ \Rightarrow \ (x,y)\in \rho$,
\item[--]  \emph{relation preserving} if $(f (x), f(y))\in \sigma
\ \Leftrightarrow \ ((x,y)\in \rho$,
\end{itemize}
\noindent for any $x,y\in S$. Clearly $f$ is relation preserving if
and only if it is isotone and reverse isotone.

A mapping $f : S \rightarrow T$ between two relational structures
$(S,\rho)$ and $(T,\sigma)$ is called an \emph{order bijection} if
it is an isotone and reverse isotone bijection.  We will write
$S\cong_o T$ if there is an order bijection $f : S \rightarrow T$
between two relational structures $(S,\rho)$ and $(T,\sigma)$.

Any mapping $f : S \rightarrow T$  between a set $S$ and relational
system $(T,\sigma)$ gives rise to the following important relation
on its domain:
$$\eta_f \stackrel{\rm def}{=}
\{(x,y)\in S \times S : (f (x), f(y))\in \sigma \}.$$

  Certain connections
between $\eta_f$ and the basic relation $\rho$ of the domain set
$S$, given in the lemma below, are, from a practical point of view
useful criteria in recognizing types of monotonicity.
\begin{lemma} \label{pl1}
 Let $(S,\rho)$ and $(T,\sigma)$  be  relational systems, and let $f : S \rightarrow
T$.  Then
 \begin{itemize}
  \item[\emph{(i)}] \quad  $f$  is isotone iff $\rho \subseteq \eta_f$;
 \item[\emph{(ii)}]  \quad   $f$  is reverse isotone iff
$\eta_f \subseteq \rho $;
 \item[\emph{(iii)}] \quad $f$  is relation preserving iff
$\eta_f = \rho $.
\end{itemize}
\end{lemma}

Let $(S,\rho)$ be a relational system and let $\varepsilon$ be an
equivalence on $S$. Let us define a binary relation
$\Theta_{\rho}=\rho/\varepsilon$ on the set $S/\varepsilon$ as
$$(x\varepsilon, y\varepsilon)\in\rho/\varepsilon \stackrel{\rm def}{\Leftrightarrow}
 (x,y)\in\rho.$$
The relational system $(S/\varepsilon; \rho/\varepsilon)$ is called
a \emph{relational factor system } of $S$ by $\varepsilon$, or a
\emph{quotient relational system} of  $(S,\rho)$ by $\varepsilon$.

\begin{lemma}\label{rel-structure}
Let $(S,\rho)$ be a reflexive, antisymmetric or transitive ordered
system and let $\varepsilon$ be an equivalence on $S$. Then
$(S/\varepsilon, \Theta_\rho)$ also a reflexive, antisymmetric or
transitive ordered system.
\end{lemma}

In what follows the following elementary result  will be handy,
 and will be used without reference.
\begin{lemma}\label{msart-14}
Let $\rho$ be a relation on a set $S$. If $\rho$ is reflexive
(symmetric, antisymmetric, transitive) then $\rho^{-1}$ is also
reflexive (symmetric, antisymmetric, transitive).
\end{lemma}
Due to Lemma~\ref{rel-structure}, for a given quasiordered set
$(S,\rho)$ and an equivalence $\varepsilon$ defined on it,
$(S/\varepsilon; \Theta_\rho)$ is a quasiordered set too. The next
lemma, which gives a way  to construct an ordered set from any given
quasiorder, is a version of the  celebrated \emph{Birkhoff's lemma},
see (Birkhoff, 1967, Lemma 1, p.  21).

\begin{lemma}\label{msart-bl}
Let $(S,\rho)$ be a quasiordered set. Then:
\begin{itemize}
\item[\emph{(i)}] \quad $\varepsilon_{\rho}= \rho\cap\rho^{-1}$ is the greatest equivalence on $S$
contained in a quasiorder $\rho$;
\item[\emph{(ii)}] \quad The relation $\hat{\Theta}_\rho=\rho/\varepsilon_{\rho}$ defined on
$S/\varepsilon_{\rho}$ by
$$(x\varepsilon_{\rho}, y\varepsilon_{\rho})\in\rho/\varepsilon_{\rho}
=\hat{\Theta}_\rho  \stackrel{\rm def}{\Leftrightarrow}
(x,y)\in\rho$$ is an ordering relation, i.e. $(S/\varepsilon_{\rho};
\hat{\Theta}_\rho)$ is an ordered set.
\item[\emph{(iii)}] \quad  $\pi_{nat}: S\rightarrow S/\varepsilon_{\rho}$ is a
surjective isotone and reverse isotone mapping between the
quasiordered set $(S,\rho)$ and the ordered set
$(S/\varepsilon_{\rho}, \hat{\Theta}_\rho)$.
\end{itemize}
\end{lemma}

The order $\hat{\Theta}_\rho\ (=\rho/\varepsilon_{\rho})$ is called
the  \emph{order defined by the quasiorder} $\rho$. Sometimes
$\hat{\Theta}_\rho$ is referred to as the \emph{kernel of the
quasiorder} $\rho$. In addition, because of the
Lemma~\ref{msart-bl}, a quasiorder is often called a
\emph{preorder}. Let us remember that a quasiorder $\alpha$ on $S$
is an order if and only if $\alpha\cap\alpha^{-1}= \triangle_S$.

So far we saw that the relation $\Theta_\rho$, defined on the
quotient set $S/\varepsilon$ by an equivalence $\varepsilon$ of a
set $S$, is defined in the terms of representatives of the
equivalent classes. In the lemma which follows it will be shown that
the definition of  $\Theta_\rho$ to a certain extend does not depend
on the representatives.
\begin{lemma}\label{msart-15}
Let $(S,\rho)$ be a quasiordered set. Then the following conditions
are equivalent:
\begin{itemize}
\item[\emph{(i)}] \quad  $(x\varepsilon, y\varepsilon)\in\Theta_{\rho}
\Leftrightarrow (\exists a\in x\varepsilon_{\rho})(\exists b\in
y\varepsilon_{\rho})\, (a,b)\in\rho$;
\item[\emph{(ii)}] \quad $(x\varepsilon_{\rho}, y\varepsilon_{\rho})\in \Theta_{\rho} \Leftrightarrow
 (x,y)\in\rho$;
\item[\emph{(iii)}] \quad $(x\varepsilon_{\rho}, y\varepsilon_{\rho})\in\Theta_{\rho}
\Leftrightarrow (\forall a\in x\varepsilon_{\rho})(\forall b\in
y\varepsilon_{\rho})\, (a,b)\in\rho$.
\end{itemize}
\end{lemma}

Let $(S,\alpha)$ be an ordered set and let $\varepsilon$ be an
equivalence defined on it. Then, in general, the relation
$\Theta_{\alpha}$ need not be an order on the factor set
$S/\varepsilon$. The following example,   (Kehayopulu and Tsingelis,  1995),
provides evidence for that.
\vspace{.5cm}

\textbf{Example 1}~~
Let $(S,\alpha)$ be an ordered set with $S=\{a,b,c,d,e\}$ and
$$\alpha = \{(a,a), (a,d), (b,b), (c,c), (c,e),(d,d), (e,e).\}$$
Let $\varepsilon$ be an equivalence relation defined on $S$ as:
$$\varepsilon = \{(a,a), (a,e), (b,b), (c,c), (c,d), (d,c), (d,d), (e,a)
(e,e)\}.$$ So $ S/\varepsilon = \{\{a,e\}, \{c,d\}, \{b\}\}$, and
let $\Theta_{\alpha}$ be a relation on $S/\varepsilon$ defined, as
usual,
$$(x\varepsilon, y\varepsilon)\in \Theta_{\alpha} \Leftrightarrow
(\exists a\in x\varepsilon)(\exists b\in y\varepsilon)\,
 (a,b)\in\alpha.$$
Then we have $(a\varepsilon, c\varepsilon)\in \Theta_{\alpha}$ since
$(a,d)\in \alpha$ as well as $(c\varepsilon,a\varepsilon)$ since
$(c,e)$ and $a\varepsilon\neq c\varepsilon$, i.e. $\Theta_{\alpha}$
is not antisymmetric. Thus $(S/\varepsilon,\Theta_{\alpha})$ is not
an ordered set. $\diamond$
\vspace{.5cm}

The \emph{First isomorphism theorem for ordered sets} follows.
\begin{theorem}\label{d-msart-11}
Let $f:S \rightarrow T$ be a mapping between ordered sets $(S,
\rho)$ and  $(T,\sigma)$. Then the following statements are true:
\begin{itemize}
\item[\emph{(i)}] \quad relation $\eta_f$ is a quasiorder on $S$.
\item[\emph{(ii)}] \quad $ker\, f= \varepsilon_f =\eta_f\cap \eta_f^{-1}$.
\item[\emph{(iii)}] \quad  The relation $\Theta_f = \eta_f/\varepsilon_f$
is an order on $S/\varepsilon_f$. The mapping $\varphi :
S/\varepsilon_f \rightarrow T$ defined by $\varphi(x(ker\,f)) =
f(x)$ is the unique isotone injective mapping between  ordered sets
$(S/\varepsilon_f, \Theta_f)$ and $(T,\sigma)$, such that
$f=\varphi\circ\pi_{nat}$ and $S/\varepsilon_f \cong_o f(S)$.
\item[\emph{(iv)}] \quad  If $f$ is surjective, then $\varphi$ is an order
bijection, i.e. we have $S/\varepsilon_f\cong_o T$.
\end{itemize}
\end{theorem}
As a consequence of Theorem~\ref{msart-12} and
Theorem~\ref{d-msart-11} we have next result.

\begin{theorem}\label{d-msart-112}
Let $f:S \rightarrow T$ be a mapping between ordered sets $(S,
\rho)$ and  $(T,\sigma)$. Let $\rho$ be a quasiorder on $S$ such
that $\rho\subseteq\eta_f$. Then the following statements are true:
\begin{itemize}
\item[\emph{(i)}] \quad $\varepsilon_{\rho} \subseteq ker\, f$.
\item[\emph{(ii)}] \quad The mapping $\varphi : S/\varepsilon_{\rho} \rightarrow T$
defined by $\varphi(x\varepsilon_{\rho}) = f(x)$ is the unique
isotone injective mapping between  ordered sets
$(S/\varepsilon_{\rho},\hat{\Theta}_\rho)$ and $(T,\sigma)$,
 such that $f=\varphi\circ\pi_{\rho}$.
\end{itemize}
\end{theorem}
To conclude:
\begin{itemize}
\item Two isomorphism theorems for sets (Theorem~\ref{msart-11}
and Theorem~\ref{msart-12}) are based on equivalences;
\item In the case of ordered sets, quasiorders play the role of
equivalences (Theorem~\ref{d-msart-11} and
Theorem~\ref{d-msart-112}).
\end{itemize}

\subsection{Semigroups within CLASS}\label{ch-msa-22}

Although frequently fragmentary, studies of semigroups carried out
at the end of 1920s  are often marked as the beginning of their
history. In the era without fast communication devices, semigroups
took their place in development of other mathematical disciplines
pretty fast. ``The analytical theory of semi-groups is a recent
addition to the ever- growing list of mathematical disciplines. It
was my good fortune to take an early interest in this discipline and
see its reach maturity. It has been a pleasant association: I hail a
semigroup when I see one and I seem to see them anywhere! Friends
have observed however that there are mathematical objects which are
not semi-groups.'' See Foreword in Hille's book (1948).

\subsubsection{Basic concepts - definitions, examples, embedding
theorems}\label{ch-msa-221}


A \emph{semigroup}  $(S, \, \cdot)$ is a nonempty set $S$  with a
binary operation $\cdot$ called  multiplication such that:
$$(x\,\cdot\, y)\,\cdot z\, = \, x\,\cdot\,(y\,\cdot\,z ),$$
for any $x,y,z \in S$.  Frequently, $xy$ is written rather than
$x\,\cdot\,y$.

A mapping $f : S \rightarrow   T$ between  two semigroups $S$ and
$T$, that preserves the operations or  is compatible with the
operations of the semigroups, i.e.
$$f(xy)= f(x)f(y)$$
is called a \emph{homomorphism}. Several types  of homomorphisms
have a specific name. A homomorphism $f$ is

$\cdot$  \emph{embedding} or \emph{monomorphism:} $f$ is  one-one;

$\cdot$  \emph{onto} or \emph{epimorphism}: $f(S)= T$;

$\cdot$  \emph{isomorphism}: $f$ is onto embedding.

\noindent If there is an isomorphism $f : S \rightarrow   T$ then we
say $S$ and $T$ are \emph{isomorphic} and denote this by $S \cong
T$.

A homomorphism $f : S \rightarrow   S$ is called an
\emph{endomorphism}, and if it is, in addition,  bijective then it
is called an \emph{automorphism}.

Some of the fundamental examples of semigroups will be presented.
Semigroups whose elements are  the binary relations or the
transformations  on a given set play a central role in applications.

Let $X$ be a set. Let $\rho, \sigma\in\mathcal{B}_X$
(Section~\ref{ch-msa-21}). Define the composition of $\rho$ and
$\sigma$ to be
$$\rho\circ\sigma =\{(x,y): (\exists z\in X) (x,z)\in\rho \
\text{and }\ (z,y)\in \sigma\}.$$
\begin{proposition}\label{ch-msa-200}
$(\mathcal{B}_X, \circ)$  is a monoid.
\end{proposition}
In what follows we shall write, as usual, $\rho\,\sigma$ instead of
$\rho\circ\sigma$.  In addition,  the ordered pair $(\mathcal{B}_X ,
\circ)$ will be denoted just as $\mathcal{B}_X$. Recall, a monoid is
a semigroup with the identity.

The semigroup $\mathcal{T}_X$ is called the \emph{full
transformation semigroup} on $X$. A subsemigroup of $\mathcal{T}_X$
is called a \emph{transformation semigroup} on $X$. The group
$\mathcal{G}(X)$ is called the \emph{symmetric group} on $X$. It is
stated that permutation groups can serve as models for all groups.
The fundamental importance of permutations group in algebra follows
from the well-known  Cayley's theorem for groups.

\begin{theorem}
If $G$ is a group then there exists an embedding $\varphi:
G\rightarrow \mathcal{G}_X$ for some $X$.
\end{theorem}

Now, The Cayley theorem  for semigroups follows.

\begin{theorem}
If $S$ is a semigroup then there exists an embedding $\varphi:
G\rightarrow \mathcal{T}_X$ for some $X$.
\end{theorem}

A homomorphism $\varphi$ from $S$ into some $\mathcal{T}_X$ is
called a \emph{representation of} $S$ \emph{by mappings}. A
representation $\varphi$ is called \emph{faithful} if it is an
embedding.

\subsubsection{General structure results -
quotient semigroups, homomorphism and isomorphism
theorems}\label{ch-msa-222}

A relation $\rho$ defined on a semigroup $S$ is called
\begin{itemize}
\item \emph{left compatible} (with  multiplication): $(x,y)\in \rho
\Rightarrow (zx,zy)\in\rho$,
\item \emph{right compatible} (with  multiplication): $(x,y)\in \rho
\Rightarrow (xz,yz)\in\rho$,
\item  \emph{compatible} (with  multiplication): $(x,y), (s,t)\in \rho
\Rightarrow (xs,yt)\in\rho$,
\end{itemize}
for any $x,y,z,s,t\in S$.

\begin{lemma}\label{ch-msart-3-l1}
Let $S$ be a semigroup and $\rho$ a quasiorder defined on it. Then,
$\rho$ is compatible if and only if it is left and right compatible.
\end{lemma}

A (left, right) compatible equivalence relation $\varepsilon$ on a
semigroup $S$ is called a (\emph{left}, \emph{right})
\emph{congruence}. The quotient set $S/\varepsilon$ is then provided
with a semigroup structure.

\begin{theorem}\label{sth21}
Let $S$ be a semigroup and $\varepsilon$ a congruence on it.Then
$S/\varepsilon$ is a semigroup with respect to the operation defined
by  $(x\varepsilon)(y\varepsilon) = (xy)\varepsilon$, and the
mapping $\pi\,:\, S\,\rightarrow\, S/\varepsilon$,
$\pi(x)=x\varepsilon$, $x\in S$, is an  epimorphism.
\end{theorem}
The semigroup $S/\varepsilon$ from the Theorem~\ref{sth21} is called
a  \emph{quotient} or \emph{factor semigroup }of $S$. The
epimorphism $\pi$ is called the \emph{natural} or \emph{canonical
epimorphism} associated with the congruence $\varepsilon$.

The theorem which follows can be considered as a consequence of the
Theorem~\ref{msart-11} and Theorem~\ref{msart-12}.

\begin{theorem}\label{sth2-n}
Let $f: S \rightarrow T$ be a homomorphism between semigroups $S$
and $T$.  Then
\begin{itemize}
\item[\emph{(i)}] \quad $ker\,f = f\circ f^{-1}$ is a congruence on $S$;
and the mapping $\varphi \,:\, S/ker\,f \, \rightarrow \, T$ defined
by $\varphi(x(ker\,f)) = f(x)$ is  an embedding such that
$f=\varphi\circ\pi$. If $f$ an epimorphism, then $\varphi$ is an
isomorphism.
\item[\emph{(ii)}] \quad  If $\varepsilon$ is a congruence on a semigroup  $S$ such
that $\varepsilon\subseteq ker\,f$, then there exists a homomorphism
of semigroups $\varphi : S/\rho  \rightarrow T$, such that
$f=\varphi\circ\pi$.
\end{itemize}
\end{theorem}


\subsubsection{Ordered semigroups -  definitions, ordered homomorphism and ordered
isomorphism theorems}\label{ch-msa-224}


Let $(S, \cdot)$ be a semigroup and $\rho$ relation on it. A triple
$(S, \cdot;\, \rho)$ is called:
\begin{itemize}
\item a \emph{quasiordered semigroup} if $\rho$ is a compatible quasiorder,
\item an \emph{ordered  semigroup} if $\rho$ is a compatible order.
\end{itemize}
Of course, a quasiorder or order defined on a semigroup $S$ need
not, in general, be compatible.

Now, Birkhoff's lemma for semigroups follows next.

\begin{lemma}\label{msart-3-bl}
Let $(S,\rho)$ be a quasiordered semigroup. Then:
\begin{itemize}
\item[\emph{(i)}] \quad $\varepsilon_{\rho}= \rho\cap\rho^{-1}$ is the greatest congruence on $S$
contained in a quasiorder $\rho$;
\item[\emph{(ii)}] \quad  $(S/\varepsilon_{\rho}; \hat{\Theta}_\rho)$ is an
ordered semigroup where
$$(x\varepsilon_{\rho}, y\varepsilon_{\rho})\in\rho/\varepsilon_{\rho}
= \hat{\Theta}_\rho  \stackrel{\rm def}{\Leftrightarrow} (x,y)\in\rho.$$
\item[\emph{(iii)}] \quad $\pi: S\rightarrow S/\varepsilon_{\rho}$ is an
 isotone and reverse isotone epimorphism  between the quasiordered
semigroup $(S,\rho)$ and the ordered semigroup
$(S/\varepsilon_{\rho}, \hat{\Theta}_\rho)$.
\end{itemize}
\end{lemma}

In addition, because of Lemma~\ref{msart-3-bl} for quasiordered
semigroup, a compatible quasiorder  is sometimes called a
\emph{half-congruence}.

The results which follow are consequences of the ones given in
Subsection~\ref{ch-msa-213} as well as in
Subsection~\ref{ch-msa-222}. The \emph{First isomorphism theorem for
ordered semigroups} follows first.
\begin{theorem}\label{d-msart-3-11}
Let $f:S \rightarrow T$ be a homomorphism between ordered semigroups
$(S, \rho)$ and  $(T,\sigma)$. Then the following statements are
true:
\begin{itemize}
\item[\emph{(i)}] \quad  $ker\, f= \varepsilon_f\subseteq \eta_f$.
\item[\emph{(ii)}] \quad  The mapping $\varphi : S/ker\,f \rightarrow T$
defined by $\varphi(x(ker\,f)) = f(x)$ is the unique isotone
embedding of the ordered  semigroup $(S/ker\,f, \Theta_f)$ into
$(T,\sigma)$, such that $f=\varphi\circ\pi$ and $S/ker\,f \cong_o
f(S)$.
\item[\emph{(iii)}] \quad  If $f$ is epimorphism, then $\varphi$ is an order
isomorphism, i.e. we have $S/ker\,f\cong_o T$.
\end{itemize}
\end{theorem}

\begin{theorem}\label{d-msart-3-111}
Let $f:S \rightarrow T$ be a homomorphism between ordered sets $(S,
\rho)$ and $(T,\sigma)$. Let $\rho$ be a compatible quasiorder on
$S$ such that $\rho\subseteq\eta_f$. Then the following statements
are true:
\begin{itemize}
\item[\emph{(i)}] \quad $\varepsilon_{\rho} \subseteq ker\, f$.
\item[\emph{(ii)}] \quad  The mapping $\varphi : S/\varepsilon_{\rho} \rightarrow T$
defined by $\varphi(x\varepsilon_{\rho}) = f(x)$ is the unique
isotone embedding between  ordered sets
$(S/\varepsilon_{\rho},\Theta_{\rho})$ and $(T,\sigma)$,
 such that $f=\varphi\circ\pi_{\rho}$.
\end{itemize}
\end{theorem}

\subsection{Applications and possible
applications}\label{ch-msa-23}

Kline (1967) stated  that ``One can look at mathematics
as a language, as a particular kind of logical structure, as a body
of knowledge about number and space, as a series of methods for
deriving conclusions, as the essence of our knowledge of the
physical world, or merely as an amusing intellectual activity [...]
Practical, scientific, philosophical, and artistic problems have
caused men to investigate mathematics.''

Within the previous subsection we  presented some results concerning
ordered sets, semigroups and ordered semigroups within classical
settings. The capability and flexibility of the just mentioned
binary structures from the point of view of modeling and
problem-solving in extremely diverse situations have been already
pointed out, and interesting new algebraic ideas arise with binary
applications  and connections to other areas of mathematics and
sciences. In what follows we are going to turn our attention to
their application within social sciences and humanities. However, we
do not pretend to mention all of the existing applications.

Order and ordered structures  enter (among others fields) into
computer science and into humanities and social sciences in many
ways and on many different levels. Order enters into the
classifications of objects in two rather different levels:
\begin{itemize}
\item[$\bullet$] classifications of certain ordered sets according to various
criteria;
\item[$\bullet$] the discipline of \emph{concept analysis} provides, on a deeper
level, a powerful technique for classifying and for analysing
complex sets of data.
\end{itemize}
Recall, formal concept analysis, \textbf{FCA}, invented by Rudolf
Wille in the early 80s,   built on the mathematical theories  of
ordered sets and lattices, is based on the mathematization of
concept and concept hierarchy. ``Hierarchies occur often both within
mathematics and the `real' world and the theory of ordered sets (and
lattices) provides a natural setting in which to discuss and analyse
them,'' as written in
(Davey and Priestley, 1990, 2002).
\textbf{FCA} has been
proven successful in a wide range of applications:  artificial
intelligence,  software engineering, chemistry, biology, psychology,
linguistics, sociology. More on applications of order theory within
social sciences can be found in Davey and Priestley's book
(1990, 2002).

For application within linguistics, see, for example, Partee, ter
Meulen and Wall's book (1990). More precisely,as
 Part C (p. 249-316) of that book leads from the notions of order and to algebraic
structures such as groups, semigroups, and monoids, and on to
lattices and Boolean and Heyting algebras, which have played a
pretty  central role in certain work in the semantics of events,
mass terms, collective vs. distributive actions, etc.

On the other hand, G. Birkhoff (1971)  wrote: ``I
do not wish to exaggerate the importance for computer science of
lattices (including Boolean algebras), or of binary groups and
fields. All of these have a quite special structure. A much more
general class of algebraic systems is provided by semigroups, which
are indeed basic for a great part of algebra''. Almost at the same
time B. M. Schein (1970) stated that  ``One meets ...
semigroups much more often than groups (and much more often than one
thinks he does), and only the existing polarization of mind on the
topic of groups and other classic algebraic structures prevents one
from seeing semigroups in various processes and phenomena of
mathematics and the universe.''

In 1991 John Paul Boyd's book  \emph{Social semigroups a unified
theory of scalling and blockmodelling  as applied
to social networks } (1991) appeared. One of the main goals of that book
was, as Boyd wrote, ``to equip the reader with powerful conceptual
and analytical tools that can be used to solve other problems in the
social sciences.'' Standard concepts of semigroup theory  acquire a
concrete meaning in terms of human kinship.

``Social semigroups?  Are you serious?'' shouted B. M. Schein (1997)  in his
review of Boyd's book,  and proceeded
``semigroups do appear very naturally in various studies known under
the common umbrella name of 'social sciences.' Yet, to see them, one
needs eyes,'' .

Recall, social networks are collections of social or interpersonal
relationships linking individuals in a social group. The partially
ordered semigroup of a network was introduced as an algebraic
construction fulfilling the following two important requirements:
\begin{itemize}
\item[$\bullet$] it should be multirelational and so encompass
different types of network relations in the description of an
individual's social environment;
\item[$\bullet$]  it should be concerned with the description of
different kinds of network paths.
\end{itemize}
More about applications of partially ordered semigroups within the
area can be found in Pattison (1993).

Abstract algebra can provide a great framework for analyzing music
and abstracting the relationships found in modern (Western) music
theory to uncover other possible music. More about the applications
of semigroups in music can be found, for example, in the work of
Bras-Amor\'{o}s (2019, 2020).

At the very end of our short journey through the applications of
semigroups within the areas of social sciences, humanities and
music, let us point out again that the list of applications given
above is far from listing  all of the existing applications of
semigroup theory. This type of research can be a topic on its own
for certain types of scientific writings.

\subsection{Notes}\label{ch-msa-24}

Results presented in this section to a great extent can be seen as a
kind of generalization of these given in Chajda (2004),
Chajda and  Sn\'{a}\v sel  (1998),
K\'{o}rtesi, Radeleczki and Szil\'{a}gyi (2005), for relational
structures of a certain type; Bloom (1976),
Sz\'{e}dli and Lenkehegyi (1983, 1983a), for the
ordered algebras;  and Kehayopulu and Tsingelis
(1995, 1995a)  for the ordered
semigroup case.

More about ordered sets and their history can be found, for example,
in Duffus and Rival (1981), Schr\"{o}der (2003).
 The standard reference for semigroup theory is
Howie  (1995). For ordered semigroup theory and on some
other algebraic structures see Blyth (2005).

\section{Within BISH}\label{ch-msa-3}

Throughout this chapter \emph{constructive mathematics} is
understood as mathematics performed in the context of intuitionistic
logic, that is, without the law of the excluded middle
(\textbf{LEM}). \textbf{LEM} can be regarded as the main source of
nonconstructivity. It was Brouwer (1975),
 who first observed that  \textbf{LEM} was extended without
justification to statements about infinite sets.
 By constructive mathematics we mean Bishop-style mathematics,
\textbf{BISH}. Several consequences of \textbf{LEM} are not accepted
in Bishop's constructivism. We will mention two such nonconstructive
principles - the ones which will be used latter.
\begin{quote}
$\bullet$  \textbf{The limited principle of omniscience},
\textbf{LPO}: \ for each binary sequence $\left( a_{n}\right)
_{n\geq1}$, either  $a_{n}=0$ for all $n\in \mathbb{N}$, or else
there exists $n$ with $a_{n}=1$.
\end{quote}%
\begin{quote}
$\bullet$  \textbf{ Markov's principle}, \textbf{MP}: \ For each
binary sequence $\left( a_{n}\right) _{n\geq1}$, if it is impossible
that $a_{n}=0$ for all $n\in \mathbb{N}$, then there exists $n$ with
$a_{n}=1$.
\end{quote}%
{\rm\textbf{LPO}} is equivalent  to the decidability of equality on
the real number line $\mathbb{R}$.
$$\forall_{x\in \mathbb{R}}\, (x=0 \vee x\neq 0).$$
Within  constructive mathematics,  a statement $P$, as in classical
mathematics, can be disproved by giving a counterexample. However,
it is also possible to give a \emph{Brouwerian counterexample} to
show that the statement is nonconstructive. A Brouwerian
counterexample to a statement $P$ is a constructive proof that $P$
implies some nonconstructive principle, such as, for example,
\textbf{LEM}, and its weaker versions \textbf{LPO}, \textbf{MP}. It
is not a counterexample in the true sense of the word - it is just
an indication that $P$ does not admit a constructive proof.

Following Troelstra and  van Dalen (1998),
\emph{constructive algebra} is more complicated than classical
algebra in various ways: algebraic structures as a rule do not carry
a decidable equality relation -- this difficulty is partly met by
the introduction of a strong inequality relation, the so-called
\emph{apartness }relation; there is (sometime) the awkward abundance
of all kinds of substructures, and hence of quotient structures.

As highlighted by M. Mitrovi\'c  {\it et al} (2019) ,
 \emph{the theory of constructive semigroup with apartness} is a new
approach to semigroup theory and not a new class of semigroups. It
presents a semigroup facet of some relatively well established
direction of constructive mathematics which, to the best of our
knowledge, has not yet been considered within the semigroup
community.

\subsection{Set with apartness}\label{ch-msa-31}

The cornerstones for \textbf{BISH} include the notion of positive
integers, sets and functions. The set $\mathbb{N}$ of positive
numbers is regarded as a basic set, and it is assumed that the
positive numbers have the usual algebraic and order properties,
including  mathematical induction.  Contrary to the classical case,
a set exists only when it is defined.

\subsubsection{Basic concepts and important examples}\label{ch-msa-311}

To define a set  $S$, we have to give a property that enables us to
construct members of $S$, and to describe the equality $=$ between
elements of $S$,  which is a matter of convention, except that it
must be an  equivalence. A set $(S, =)$ is an \emph{inhabited} set
if we can construct an element of $S$. The distinction between the
notions of a nonempty set and an inhabited set is a key in
constructive set theories. The notion of equality of different sets
is not defined. The only way in which elements of two different sets
can be regarded as equal is by requiring them to be subsets of a
third set.

A property $P$, which is applicable to the elements of a set $S$,
 determines a subset of $S$ denoted by $\{x\in S : P(x)\}$.
 Furthermore, we will be interested only in properties $P(x)$
which are \emph{extensional} in the sense that for all $x_1,x_2\in
S$ with $x_1 = x_2$, $P(x_1)$ and $P(x_2)$ are equivalent.
Informally, it means that ``it does not depend on the particular
description by which $x$ is given to us''.

Let $(S,=)$ be an \emph{inhabited} set. By an
\textbf{\emph{apartness}} on $S$  we mean a binary relation $\# $ on
$S$ which satisfies the axioms of  irreflexivity, symmetry and
cotransitivity:
\begin{itemize}
\item[(Ap1)] \qquad $\neg(x \# x)$

\item[(Ap2)] \qquad $x \# y \ \Rightarrow\ y \# x$,

\item[(Ap3)] \qquad $x\# z\ \Rightarrow\ \forall _y\,(x\# y\,\vee\,y\# z).$
\end{itemize}
If $x\# y$, then $x$ and $y$ are different, or distinct.
 Roughly speaking, $x= y$ means that we have a proof that $x$ equals $y$ while
$x\# y$ means that we have a proof that $x$ and $y$ are different.
Therefore, the negation of $x= y$ does not necessarily imply that
$x\# y$ and vice versa: given $x$ and $y$, we may have neither a
proof that $x= y$ nor a proof that $x\# y$.

The apartness on a set $S$ is \emph{tight} if
\begin{itemize}
\item[(Ap4)] \qquad $\neg(x \# y) \ \Rightarrow\ x=y$.
\end{itemize}
If apartness is tight, then  $\neg(x \# y) \ \Leftrightarrow\ x=y$.

By extensionality, we have

\begin{itemize}
\item[(Ap5)] \qquad $x \# y \,\wedge\,y = z \ \Rightarrow\, x \# z$,
\end{itemize}
\noindent the equivalent form  of which is

\begin{itemize}
\item[(Ap5')] \qquad $x \# y \,\wedge\,x = x' \,\wedge\, y=y'\ \Rightarrow\,
x' \# y'$.
\end{itemize}

A \emph{set with apartness} $(S,=,\#)$ is the starting point for our
 considerations, and will be simply denoted by $S$.
The existence of an apartness relation on a structure often  gives
rise to an  apartness relation on another structure. For example,
given two sets with apartness $(S,= _S,\#_S)$ and $(T,=_T,\#_T)$, it
is permissible to construct the  set of mappings between them. A
\emph{mapping} $f:S\rightarrow T$ is an algorithm which produces an
element $f(x)$ of $T$ when applied to an element $x$ of $S$, which
is extensional, that is
\[
\forall_{x,y\in S}\,(x=_{S}y\  \Rightarrow \ f(x)=_{T}f(y)).\mathbf{\ }%
\]

An important property applicable to mapping $f$ is that of strong
extensionality. Namely, a  mapping $f:S\rightarrow T$ is a
\emph{strongly extensional} mapping, or, for short, an
\emph{se-mapping}, if
\[
\forall_{x,y\in S}\, (f(x)\#_{T} f(y)\ \Rightarrow\ x\#_{S} y).
\]
Let us remember that strong extensionality of all mappings from
$\mathbb{R}$ to $\mathbb{R}$ implies the Markov principle,
\textbf{MP}.

An se-mapping $f$ is:

- \, an \emph{se-surjection} if it is surjective;

- \, an  \emph{se-injection} if it is injective;

- \, an  \emph{se-bijection} if it is bijective.

 Furthermore, $f$ is

 -  \emph{apartness injective}, shortly \emph{a-injective}:  $\forall_{x,y\in S}\, (x\#_{S}\, y\  \Rightarrow \
f(x)\#_{T} f(y))$;

-  \emph{apartness bijective}:  a-injective, se-bijective.

Given a set $X$ with apartness it is permissible to construct the
set of all se-mappings of $X$ into itself which inherits apartness
from $X$.

\begin{theorem}  \label{ch-msa-thm1.2}
 Let $(X, =, \#)$ be a set with apartness.
 If $\mathcal{T}_S^{se}=M(X,X)$ is the set of all se-mappings from $X$ to $X$,
 then $\pmb{(\mathcal{T}_S^{se}, =, \#)}$ with
\begin{eqnarray*}
f=g \ &\Leftrightarrow& \  \forall _{x\in X} \, (f(x)=g(x)) \\
f\# g \  &\Leftrightarrow&\ \exists _{x\in X} \, (f(x)\#  g(x))
\end{eqnarray*}
is a set with apartness too.
\end{theorem}

Given two sets with apartness $S$ and $T$ it is permissible to
construct the set of ordered pairs  $(S\times T, =, \#)$ of these
sets defining apartness by $$(s,t)\,\#\, (u,v) \enspace
\stackrel{\rm def}{\Leftrightarrow}\enspace s \,\#_S\, u\,\vee \,
t\, \#_T\, v.$$

\subsubsection{Distinguishing subsets}\label{ch-msa-312}

The presence of apartness implies the appearance of different types
of substructures connected to it. Following Bridges and
V\^{\i}\c{t}\={a} (2011), we define the relation
$\bowtie$ between an element $x\in S$ and a subset $Y$ of $S$  by
$$x\bowtie Y \ \stackrel{\rm def}{\Leftrightarrow} \ \forall _{y\in Y}  (x\# y).$$
A subset $Y$  of $S$ has two natural complementary subsets:
\emph{the logical complement} of $Y$
$${\neg Y} \stackrel{\rm def}{=} \{x\in S : x\notin Y\}, $$ and \emph{the apartness complement}
 or, shortly,
the \emph{a-complement} of $Y$ $${\sim Y} \,  \stackrel{\rm def}{=}
\, \{x \in S\, : \, x \bowtie Y\}.$$ Denote by $\widetilde{x}$ the
a-complement of the singleton $\{ x\}$. Then it can be easily shown
that $x\in \sim Y$ if and only if $Y\subseteq \widetilde{x}.$

If the apartness is not tight we can find subsets $Y$ with ${\sim Y}
\subset {\neg Y}$ as in the following example.
\vspace{.5cm}

\textbf{Example 2}~~
Let  $S=\{a,b,c\}$ be a set with apartness defined by \\
$\{(a,c),(c,a),(b,c),(c,b)\}$ and let $Y=\{ a\}$. Then the
a-complement $\sim\!\,Y = \{ c\}$ is a proper subset of its logical
complement $\neg Y=\{ b,c\}$. $\diamond$
\vspace{.5cm}

The complements  are used for the classification of subsets of a
given set. A subset  $Y$  of $S$ is
\begin{itemize}
  \item   \emph{ a detachable} subset  in $S$  or, in short, a
\emph{d-subset} in $S$ if
$$
\forall _{ x\in S} \, (x \in Y \vee x \in {\neg} Y);
$$
\item \emph{a strongly detachable} subset of $S$, shortly \emph{an
sd-subset of} $S$, if
$$
\forall _{ x\in S} \, (x \in Y \vee x \in {\sim} Y),
$$
\item \emph{a quasi-detachable} subset of $S$, shortly \emph{a
qd-subset of} $S$, if
$$
\forall _{ x\in S} \,\forall _{ y\in Y} \, (x \in Y \vee x \# y).
$$
\end{itemize}

 A  description of the relationships between those
subsets of set with apartness, which, in turn, justifies the
\emph{constructive order theory} for sets and semigroups with
apartness we develop, is given in the next theorem.

\begin{theorem}\label{senegcom}
Let $Y$ be a subset of $S$. Then:
\begin{itemize}
\item[\emph{(i)}] \quad  Any sd-subset is a qd-subset of $S$. The converse implication entails  {\rm{\textbf{LPO}}}.
\item[\emph{(ii)}] \quad  Any qd-subset $Y$ of $S$ satisfies ${\sim Y}={\neg
Y}$.
\item[\emph{(iii)}]  \quad If any qd-subset is a d-subset, then
{\rm{\textbf{LPO}}} holds.
\item[\emph{(iv)}] \quad If any d-subset is a qd-subset, then {\rm{\textbf{MP}}} holds.
\item[\emph{(v)}] \quad  Any sd-subset is a d-subset of $S$. The converse implication
entails {\rm{\textbf{MP}}}.
\item[\emph{(vi)}] \quad  If any   subset of a set with apartness $S$ is a qd-subset, then
{\rm{\textbf{LPO}}} holds.
\end{itemize}
\end{theorem}

For all subsets $Y$ of the set with apartness $S$ for  which two
distinguished complements coincide, we will adopt the following
notation:
$$\pmb{Y^c} = \,\sim Y = \neg Y.$$

\subsubsection{Binary relations}\label{ch-msa-313}

Let $(S\times S,=,\#)$ be a set with apartness. An inhabited subset
of $S\times S$, or, equivalently, a property applicable to the
elements of $S\times S$, is called a \emph{binary relation} on $S$.
Let $\alpha$ be a relation on $S$.
 Then $$(a,b) \bowtie  \,\alpha \  \Leftrightarrow \ \forall
_{(x,y)\in\alpha }\ ((a,b)\, \#\, (x,y)),$$ \noindent for any
$(a,b)\in S\times S$. The apartness complement of $\alpha$ is the
relation
$$\mathbf{\sim \alpha } \, =\, \{(x,y)
\in S\times S : (x,y) \bowtie \alpha\}.$$ In general, we have
$\sim\alpha \subseteq \neg \alpha$, which is shown by the following
example.
\vspace{.5cm}

\textbf{Example 3}~~
Let  $S=\{a,b,c\}$ be a set with apartness defined by \\
$\{(a,c),(c,a),(b,c),(c,b)\}$. Let  $\alpha =\{(a,c),(c,a)\}$ be a
  relation on $S$. Its a-complement
$$\sim\!\,\alpha = \{(a,a),(b,b),(c,c),(a,b), (b,a)\}$$
is a proper subset of its logical complement $\neg\alpha$.
$\diamond$
\vspace{.5cm}

The subset $\bigtriangledown_S=\{(x,y)\in S\times S: \neg(x= y)\}$
of $S\times S$ is called the {\em co-diagonal} of $S$.

Let $\alpha$ and $\beta$ be relations on $S$. Then, the following
operations can be defined:
\begin{itemize}
\item[-] the \emph{composition} or \emph{product} of $\alpha$ and $\beta$:
 $$\alpha\,\circ\,\beta = \left\{(x,z)\in S\times S : \exists_{y\in S}\,
   ((x,y)\in\alpha  \wedge (y,z)\in\beta) \right\};$$
\item[-] the \emph{co-composition} or \emph{co-product} of $\alpha$ and $\beta$:
 $$\alpha\,\ast\,\beta = \left\{(x,z)\in S\times S : \exists_{y\in S}\,
   ((x,y)\in\alpha  \vee (y,z)\in\beta) \right\};$$
\item[-]
$\alpha$ is \emph{associated} with $\beta$: $$\alpha \looparrowleft
\beta \stackrel{\rm def}{\Leftrightarrow}
 \forall _{x, y,z\in S} \ ((x,y)\in \alpha \wedge (y,z)\in \beta \
\Rightarrow \  (x,z)\in \alpha ).$$
\end{itemize}
In what follows, the next lemma will be of practical use.
\begin{lemma} \label{ch-msa-*properties}
  Let $\alpha$ and $\beta$ be relations on $S$ such that
 $\alpha\subseteq\gamma$ and   $\beta\subseteq\delta$. Then
  $\alpha\ast\beta\subseteq\gamma\ast\delta$.
\end{lemma}

The relation $\alpha$ defined on a set with apartness $S$ is
\begin{itemize}
\item[$\bullet$]  strongly irreflexive if  $(x,y)\in \alpha \ \Rightarrow\ x\#y$;
\item[$\bullet$] co-transitive if $(x,y)\in\alpha \ \Rightarrow\ \forall _{z\in S}\,((x,z)\in\alpha
\,\vee\,(z,y)\in\alpha)$;
\item[$\bullet$]  \emph{co-antisymmetric} if $\forall_{ x,y\in S}\, (x\# y\Longrightarrow (x,y)\in\alpha \vee
(x,y)\in\alpha^{-1})$.
\end{itemize}
An alternatively way of defining properties given above is:
\begin{itemize}
\item[$\bullet$] \emph{irreflexive} if $\alpha\subseteq \bigtriangledown$

\item[$\bullet$]  \emph{strongly irreflexive} if  $ \alpha\subseteq\#$;

\item[$\bullet$] \emph{co-transitive} if  $\alpha\subseteq \alpha\ast \alpha$;

\item[$\bullet$] \emph{co-antisymmetric} if \ $\#\,\subseteq\, (\alpha\cup\alpha^{-1})$.
\end{itemize}
It is easy to check that a strongly irreflexive relation is also
irreflexive. For a tight apartness, the two notions of irreflexivity
are classically equivalent but not so constructively. More
precisely, if each irreflexive relation were strongly irreflexive
then\textbf{ MP} would hold.

\begin{lemma} \label{ch-msa-pl2}
Let $\alpha$ be a  relation on $S$. Then:
\begin{itemize}
 \item[\emph{(i)}] \quad  If $\alpha$ (respectively, $\beta$) is strongly irreflexive then $\alpha\ast\beta\subseteq\beta$
    (respectively, $\alpha\ast\beta\subseteq\alpha$);
\item[\emph{(ii)}] \quad $\alpha$ is strongly irreflexive then $\alpha \ast\alpha$ is strongly irreflexive.
\end{itemize}
\end{lemma}

In what follows the following result will be in handy, and will be
used without further reference.
\begin{lemma} \label{ch-msa-sepl2}
If $\alpha$  is strongly irreflexive or co-transitive then
$\alpha^{-1}$ is strongly irreflexive or co-transitive too.
\end{lemma}

Now, we turn our attention to certain properties of the a-complement
of given types of relations.
\begin{lemma} \label{ch-msa-reflexiveconsistent}
 Let $\alpha$ be a relation on $S$. Then:
 \begin{itemize}
 \item[\emph{(i)}] \quad $\alpha$ is strongly irreflexive if and only if $\sim\!\,\alpha$ is reflexive;
 \item[\emph{(ii)}] \quad  If $\alpha$ is symmetric then $\sim\!\,\alpha$ is
 symmetric;
 \item[\emph{(iii)}] \quad If $\alpha$ is co-transitive then $\sim\!\,\alpha$
 is transitive;
\item[\emph{(iv)}] \quad  If $\alpha$ is co-antisymmetric and the apartness on
$S$ is tight, then $\sim\!\,\alpha$ is antisymmetric.
\end{itemize}
\end{lemma}

The relation $\alpha$ defined on a set $S$ is
\begin{itemize}
\item[$\bullet$] \emph{co-quasiorder } if it is strongly irreflexive and
co-transitive.
\end{itemize}
The co-quasiorder is one of the main building blocks for the
co-order theory of set with (non-tight) apartness that we develop.
We can use them to define the following important types of
relations.
\begin{itemize}
\item[$\bullet$] \emph{co-equivalence}: a symmetric co-quasiorder;
\item[$\bullet$]  \emph{co-order}: a co-antisymmetric co-quasiorder.
\end{itemize}

As in Example~3 the a-complement of a relation can be a
proper subset of its logical complement. If the relation in question
is a co-quasiorder, then we have the following important properties.

\begin{proposition} \label{ch-msa-senegcompl}
 Let  $\tau$ be a co-quasiorder on $S$. Then:
\begin{itemize}
\item[\emph{(i)}] \quad   $\tau$ is a qd-subset of $S\times S$;
\item[\emph{(ii)}] \quad  $\sim\!\,\tau = \neg\,\tau (= \tau^c)$;
\item[\emph{(iii)}] \quad  $\tau^c$ is a quasiorder on $S$.
\end{itemize}
\end{proposition}

\begin{lemma} \label{ch-msa-cun}
If $\tau$  and $\sigma$ are co-quasiorders on a set $S$, then  $\tau
\cup\sigma$ is a co-quasiorder too.
\end{lemma}

\subsubsection{Apartness isomorphism theorems  for sets with apartness}\label{ch-msa-314}

A quotient structure does not have, in general, a natural apartness
relation. For most purposes, we overcome this problem using a
\textbf{\emph{co-equivalence}}--symmetric co-quasiorder--instead of
an equivalence. Existing properties of a co-equivalence guarantee
that its a-complement  is an equivalence and that the quotient set
of that equivalence will inherit an apartness. The following notion
will be necessary.  For any two relations  $\alpha$ and $\beta$ on
$S$ we can   say that $\mathbf{\alpha}$ \emph{defines an apartness
on} $\mathbf{S/\beta}$  if we have
\begin{eqnarray*}
\text{(Ap\ 6)} \quad \quad  \ \  x\beta\,\# \,y\beta \ \stackrel{\rm
def}{\Leftrightarrow}\ (x,y)\in\alpha.
\end{eqnarray*}
\begin{lemma}\label{ch-msa-l2b}
If $\alpha$ is a co-quasiorder and $\beta$  an equivalence on a set
$S$, then \rm{(Ap 6)} implies
\begin{eqnarray*}
\rm{(Ap \ 6')} \quad \quad  \ \
((x,a)\in\beta\,\wedge\,(y,b)\in\beta) \ \Rightarrow \ (
(x,y)\in\alpha \,\Leftrightarrow\, (a,b)\in\alpha).
\end{eqnarray*}
\end{lemma}

\begin{theorem}\label{ch-msa-t-coe}
Let  $\kappa$ be a co-equivalence on $S$. Then
\begin{itemize}
\item[\emph{(i)}] \quad   the relation $\kappa^c$
is an equivalence on $S$ such that $\kappa \looparrowleft \kappa^c$;
\item[\emph{(ii)}] \quad $(S/\kappa^c, =, \#)$ is a set with apartness
where
\begin{eqnarray*}
a\kappa^c\,=b\kappa^c\  & \Leftrightarrow\ (a,b)\bowtie \kappa\\
a\kappa^c\,\#\,b\kappa^c\  & \Leftrightarrow\ (a,b)\in\kappa;
\end{eqnarray*}
\item[\emph{(iii)}] \quad The quotient mapping $\pi:S\rightarrow
S/\kappa^c$, defined by $\pi(x)= x\kappa^c$, is an se-surjection.
\end{itemize}
\end{theorem}

Let $f : S\to T$ be an se-mapping between sets with apartness. Then
the relation%
\[
\mathrm{coker}\,f \stackrel{\rm def}{=} \{(x,y)\in S\times S:f(x)\#
f(y)\}
\]
defined on $S$ is called the  \emph{co-kernel }of $f$.

Now, \emph{the First apartness isomorphism theorem } for sets with
apartness follows.
\begin{theorem} \label{ch-msa-cmrth2}
Let $f : S\to T$ be an se-mapping between sets with apartness. Then
\begin{itemize}
\item[\emph{(i)}] \quad   $\mathrm{coker}\,f$ is a co-equivalence on $S$;
\item[\emph{(ii)}] \quad  $\mathrm{coker}\,f \looparrowleft
\mathrm{ker}\,f$ and $\mathrm{ker}\,f\subseteq
(\mathrm{coker}\,f)^c$;
\item[\emph{(iii)}] \quad  $(S/\ker\,f,=,\#)$ is a set with apartness,
where
\begin{align*}
a(\ker\,f)\,=\,b(\ker\,f)\   & \Leftrightarrow\ (a,b)\in\ker\,f\\
a(\ker\,f)\,\#\, b(\ker\,f)\   & \Leftrightarrow\
(a,b)\in\mathrm{coker}\,f;
\end{align*}
\item[\emph{(iv)}] \quad  the mapping $\varphi:S/\mathrm{ker}\,f
\to T$, defined by $\varphi(x(\mathrm{ker}\,f))=f(x)$, is an
a-injective se-injection such that $f=\varphi\pi$;
\item[\emph{(v)}] \quad if $f$ maps $S$ onto $T$, then
$\varphi$ is an apartness bijection.
\end{itemize}
\end{theorem}

The next theorem concerns a more general situation.
\begin{theorem}\label{ch-msa-BHM-1}
Let $S$ be a set with apartness. Then:
\begin{itemize}
\item[\emph{(i)}] \quad Let $\varepsilon$ be an equivalence, and $\kappa$ a co-equivalence on $S$. Then, $\kappa$
  defines an apartness  on the factor set $S/\varepsilon$ if and only if $\varepsilon\cap\kappa =
  \emptyset$.
\item[\emph{(ii)}] \quad  The quotient mapping $\pi : S \rightarrow S/\varepsilon$, defined by
$\pi (x)  = x\varepsilon$, is an  se-surjection.
\end{itemize}
\end{theorem}

Now,  \emph{the Second apartness isomorphism theorem},  a
generalised version of Theorem~\ref{ch-msa-cmrth2}, for sets with
apartness follows.
\begin{theorem} \label{ch-msa-basicfactor}
  Let $f:S\to T$ be a mapping between sets with apartness,  and let
  $\kappa$ be a co-equivalence on $S$ such that $\kappa\cap \mathrm{ker}\,f=\emptyset$. Then:
  \begin{itemize}
  \item[\emph{(i)}] \quad  $\kappa$ defines an apartness on the factor set $S/\mathrm{ker}\,f$;
  \item[\emph{(ii)}] \quad The projection  $\pi: S\to S/\mathrm{ker}\,f$ defined by $\pi (x)= x(\mathrm{ker}\,f)$ is an
  se-surjection;
  \item[\emph{(iii)}] \quad The mapping  $\varphi:S/\mathrm{ker}\,f\to T$, given by
  $\varphi(x(\mathrm{ker}\,f))=f(x)$,is an a-injective se-injection such that  $f=\varphi\,\pi$;
  \item[\emph{(iv)}] \quad  $\varphi$ is an se-mapping if and only if $\mathrm{coker}\,f\subseteq\kappa$;
\item[\emph{(v)}] \quad  $\varphi$ is a-injective if and only if
$\kappa\subseteq \mathrm{coker}\,f$;
\item[\emph{(vi)}] \quad  If $\varphi$ is an se-mapping, then $f$ is an
  se-mapping too.
  \end{itemize}
\end{theorem}

\subsubsection{Co-ordered sets with apartness}\label{ch-msa-315}

A \emph{constructive version of Birkhoff's lemma}
(Lemma~\ref{msart-bl}) follows.

\begin{lemma}\label{ch-msa-cbl}
Let $(S,=,\#;\, \tau)$ be a co-quasiordered set with apartness. Then
 \begin{itemize}
\item[\emph{(i)}] \quad  $\kappa_{\tau} = \tau\cup \tau^{-1}$ is a
  co-equivalence on $S$;
\item[\emph{(ii)}]\quad  $(S/\kappa_{\tau}^c,=,\#)$ is a set with apartness
where
\begin{eqnarray*}
x\kappa^c_\tau\,=y\kappa^c_\tau\  & \Leftrightarrow\ (x,y)\bowtie \kappa_\tau\\
x\kappa^c_\tau\,\#\,y\kappa^c_\tau\  & \Leftrightarrow\
(x,y)\in\kappa_\tau;
\end{eqnarray*}
\item[\emph{(iii)}] \quad $(S/\kappa^c_{\tau},=,\#;\, \Upsilon_\tau)$ is a
co-ordered set with apartness where
$$(x\kappa^c_\tau, y\kappa^c_\tau)\in\Upsilon_\tau  \stackrel{\rm
def}{\Leftrightarrow} (x,y)\in\tau;$$
\item[\emph{(iv)}] \quad The quotient mapping $\pi_\tau: S\rightarrow
S/\kappa^c_\tau$, defined by $\pi_\tau(x)= x\kappa^c_\tau$, is an
isotone and reverse isotone se-surjection.
\end{itemize}
 \end{lemma}

Any mapping between a set with apartness $(S,=,\#)$ and a
co-quasiordered set with apartness $(T,=,\#\; \sigma)$ defines the
relation $\mu_f$ on $S$ by
$$(x,y)\in\mu_f \stackrel{\rm
def}{\Leftrightarrow} (f(x),f(y))\in\sigma.$$

\begin{theorem}\label{ch-msa-cp11}
Let $(S,=,\#)$ be a set with apartness, and let $(T,=,\#;\, \sigma)$
be a co-quasiordered set with apartness. If $f:S\rightarrow T$ is an
se-mapping, then
 \begin{itemize}
\item[\emph{(i)}] \quad  $\mu_f$ is a co-quasiorder on $S$;
\item[\emph{(ii)}] \quad  $\kappa_f=\mu_f\cup\mu_f^{-1}$ is a co-equivalence on on $S$ such that
$\kappa_f\subseteq \mathrm{coker}\,f$;
\item[\emph{(iii)}] \quad If $\sigma$ is a co-order on $T$, then
$\kappa_f=\mathrm{coker}\,f$;
\item[\emph{(iv)}] \quad $(S/\kappa^c_f, =, \#; \Upsilon_f)$ where
$$(x\kappa^c_f, y\kappa^c_f)\in\Upsilon_f  \stackrel{\rm
def}{\Leftrightarrow} (x,y)\in\mu_f$$ is a co-ordered set with
apartness;
\item[\emph{(v)}] \quad The mapping $\psi: S/\kappa^c_f \rightarrow S/\kappa^c_\sigma$ between
co-ordered sets with apartness $(S/\kappa^c_f, =, \#;\Upsilon_f)$
and $(S/\kappa^c_\sigma , =, \#;\Upsilon_\sigma)$, defined by
$\psi(x\kappa_f)= (f(x))\kappa_\sigma$, is an isotone and reverse
isotone a-injective se-mapping such that $\psi\pi_f=\pi_\sigma f$.
\end{itemize}
 \end{theorem}

\begin{theorem}\label{ch-msa-cqo12}
Let $(S,=,\#,\, \tau)$ be a co-quasiordered set with apartness, and
let $(T,=,\#;\, \sigma)$ be a co-quasiordered set with tight
apartness. If $f:S\rightarrow T$ is an se-mapping such that
$\mu_f\subseteq\tau$, then the mapping $\varphi: S/\kappa_\tau^c
\rightarrow T$, defined by $\varphi(x\kappa_\tau^c)=f(x)$, is an
se-mapping such that $\varphi\pi_\tau = f$.
\end{theorem}

\subsection{Semigroups with apartness}\label{ch-msa-32}

During the implementation of the \textbf{FTA} Project, Geuvers  {\it et
al} (2002, 2002a), the notion of
commutative constructive semigroups with tight apartness appeared.
We put noncommutative constructive semigroups with ``ordinary''
apartness in the center of our study, proving first, of course, that
such semigroups do exist. The initial step towards grounding the
theory is done by our contributing papers, Mitrovi\'c and
co-authors:
(Crvenkovi\'c  {\it et al}, 2013), (Crvenkovi\'c  {\it et al}, 2016),
(Mitrovi\'c   {\it et al}, 2021),  (Mitrovi\'c and Silvestrov, 2020), (Mitrovi\'c   {\it et al}, 2019).
The theory of semigroup with apartness is a new approach to semigroup
theory and not a new class of semigroups. It presents a semigroup
facet of some relatively well established direction of constructive
mathematics which, to the best of our knowledge, has not yet been
considered within the semigroup community.

\subsubsection{Basic concepts - definitions, examples,
se-embeddings}\label{ch-msa-321}

A semigroup with apartness satisfies a number of extra conditions,
firstly the well known axioms of apartness, and secondly the
semigroup operation has to be strongly extensional.

Given a set with apartness $(S,=,\#)$, the tuple $(S, =, \#, \,
\cdot )$ is a \emph{semigroup with apartness} if the binary
operation $\cdot$ is associative
\begin{itemize}
\item[(A)] \quad  $ \enspace\quad \forall _{a,b,c \in S} \ [(a\cdot b)\cdot c\, =\, a\cdot (b\cdot c)]$,
\end{itemize}
and strongly extensional
\begin{itemize}
\item[(S)] \quad  $ \enspace\quad \forall _{a,b,x,y \in S}\ (a\cdot x \#\, b\cdot y
\Rightarrow(a \#\, b  \,\vee \,  x \#\, y))$.
\end{itemize}

\noindent As usual, we are going to write $ab$ instead of $a\cdot
b$. Example~4 provides a concrete instance of a
semigroup with apartness.
\vspace{.5cm}

\textbf{Example 4}~~
Let $S = \{a,b,c,d,e\}$ be a  set with diagonal $\triangle_S$
 as the equality relation. If we denote  by
 $K=\triangle_S \cup \{(a,b),(b,a)\}$, then we can define an apartness $\#$ on $S$
 to be  $(S\times S) \setminus  K$. Thus, $(S, =, \#)$ is a set with apartness.
If we define  multiplication on the set $S$ as
\begin{center}
\begin{tabular}{c|c c c c c } 
  \(\cdot\) & a & b & c & d & e \\ \hline
  g & a & b & c & d & e \\ \hline
  a & b & b & d & d & d \\
  b & b & b & d & d & d \\
  c & d & d & c & d & c \\
  d & d & d & d & d & d \\
  e & d & d & c & d & c \\
\end{tabular}
\end{center}
then $(S, =, \#; \cdot)$ is a semigroup with apartness. $\diamond$
\vspace{.5cm}

An important example of a class of semigroup with apartness arises
in the following way. For a given set with apartness $X$ we can
construct a semigroup with apartness $\mathcal{T}_X^{se}$ as  given
in the next theorem.

\begin{theorem}\label{ch-msa-evid}
Let $X$ be a set with apartness. Let $\mathcal{T}_X^{se}$ be the set
of all se-functions from $X$ to $X$ with
$$f=\, g \ \Leftrightarrow
\  \forall _{x\in S} \, (f(x)=\, g(x))$$ and apartness $$f\#\, g \
\Leftrightarrow\ \exists _{x\in A} \, (f(x)\#\,  g(x)).$$ Then
$(\mathcal{T}_X^{se}, =, \# , \cdot )$ is a semigroup with apartness
with respect to the binary operation of composition of functions.
\end{theorem}

Until the end of of this paper, we adopt the convention that
\emph{semigroup} means \emph{semigroup with apartness}.
 Apartness  from  Theorem~\ref{ch-msa-evid} does not have to be
tight. The following example shows that we cannot even prove
constructively that the apartness on every \emph{finite} semigroup
is tight.

\textbf{Example 5}~~
Let $X=\left\{  0,1,2\right\}  $ with the usual equality
relation, that is, the diagonal $\Delta_{X}$ of $X\times X$. Let%
\[
K=\Delta_{S}\cup\{(1,2),(2,1)\},
\]
and define an apartness $\#$ on $X$ to be $(X\times X)\setminus K$.
Then, by the Theorem~\ref{ch-msa-evid},  $\mathcal{T}_X^{se}$
becomes a semigroup with apartness. Define mappings $f,g:X\rightarrow X$ by%
\begin{align*}
f(0)  & =1,\ f(1)=1,~f(2)=2,\\
g(0)  & =2,\ g(1)=1,~g(2)=2.
\end{align*}
In view of our definition of the apartness on $X$, there is no
element $x\in X$  with $f(x)\# g(x)$; so, in particular, $f$ and $g$
are se-functions. However, if $f=g$, then $1=2$, which, by our
definition of the equality on $S$, is not the case. Hence the
apartness on $S$ is not tight. $\diamond$

Let $S$ and $T$ be semigroups with apartness. An homomorphism $f : S
\rightarrow T$ is

- an \emph{se-embedding} if it is an se-injection;

- an \emph{apartness embedding} if it is an a-injective
se-embedding;

- an \emph{apartness isomorphism} if it is an apartness bijection.

 As a consequence of Theorem~\ref{ch-msa-evid}, we can formulate \emph{the constructive
Cayley's theorem for semigroups with apartness} as follows.

\begin{theorem}\label{ch-msa-Cayley}
Every semigroup with apartness se-embeds into the semigroup of all
strongly extensional self-maps on a set.
\end{theorem}

\subsubsection{Co-quasiorders defined on a semigroup}\label{ch-msa-322}

A relation $\tau$ defined on a semigroup $S$ with apartness is
called
\begin{itemize}
\item \emph{left co-compatible}: $(zx,zy)\in \tau
\Rightarrow (x,y)\in\tau$,
\item \emph{right co-compatible}: $(xz,yz)\in \tau
\Rightarrow (x,y)\in\tau$,
\item  \emph{co-compatible}: $(xz,yt)\in \tau
\Rightarrow (x,y)\in\tau \vee (z,t)\in\tau$,
\end{itemize}
for any $x,y,z,t\in S$.

Let us proceed with an example of a co-quasiorder defined on a
semigroup with apartness $S$.
\vspace{.5cm}

\textbf{Example 6}~~
Let $S$ be a semigroup with apartness as given in the
Example~4. The  relation $\tau $, defined by
$$\tau  = \{(c,a),(c,b),(d,a),(d,b),(d,c),(e,a),(e,b),(e,c), (e,d)\},$$
is a co-quasiorder  on $S$. $\diamond$
\vspace{.5cm}

The lemma which follows will be used without reference.
\begin{lemma}\label{ch-msa-lrc}
Let $\tau$ be a co-quasiorder on a semigroup with apartness $S$.
Then, $\tau$ is co-compatible if and only if $\tau$ is a left and a
right co-compatible.
\end{lemma}

\subsubsection{Apartness isomorphism theorems for semigroups with apartness}\label{ch-msa-323}

Let us remember that in \textbf{CLASS} the compatibility property is
an important condition for providing the semigroup structure on
quotient sets. Now we are looking for the tools for introducing an
apartness relation on a factor semigroup. Our starting point is the
results from   Subsection~\ref{ch-msa-314}, as well as the next
definition.

A co-equivalence $\pmb{\kappa}$ is a \emph{co-congruence} if it is
\emph{co-compatible}

\begin{theorem}\label{ch-msa-sap2}
Let $S$ be a semigroup with apartness, and let $\kappa$ be a
co-congruence on $S$.\textbf{ }Define
\begin{align*}
a\kappa^c\,=b\kappa^c\   & \Leftrightarrow\ (a,b)\bowtie
\kappa\text{,}\\
a\kappa^c\,\#\,b\kappa^c\,  & \Leftrightarrow\ (a,b)\in
\kappa\text{,}\\
a\kappa^c\,b\kappa^c\   & =\ (ab)\kappa^c\text{.}
\end{align*}
Then $(S/\kappa^c,=,\#,\cdot\,)$ is a semigroup with apartness.
Moreover, the quotient mapping $\pi:S\rightarrow S/\kappa^c$,
defined by $\pi(x)=x\kappa^c$, is an se-epimorphism.
\end{theorem}

The \emph{First apartness isomorphism theorem for semigroups with
apartness} follows.
\begin{theorem} \label{ch-msa-cmrth22}
Let $f : S\to T$ be an se-homomorphism between sets with apartness.
Then
\begin{itemize}
\item[\emph{(i)}] \quad  $\mathrm{coker}\,f$ is a co-congruence on $S$;
\item[\emph{(ii)}] \quad  $\mathrm{coker}\,f \looparrowleft
\mathrm{ker}\,f$ and $\mathrm{ker}\,f\subseteq
(\mathrm{coker}\,f)^c$;
\item[\emph{(iii)}] \quad $(S/\ker\,f,=,\#; \cdot)$ is a semigroup with apartness,
where
\begin{align*}
a(\ker\,f)\,=\,b(\ker\,f)\   & \Leftrightarrow\ (a,b)\in\ker\,f\\
a(\ker\,f)\,\# \,b(\ker\,f)\   & \Leftrightarrow\
(a,b)\in\mathrm{coker}\,f;
\end{align*}
\item[\emph{(iv)}]  \quad The mapping $\varphi:S/\mathrm{ker}\,f
\to T$, defined by $\varphi(x(\mathrm{ker}\,f))=f(x)$, is an
apartness embedding such that $f=\varphi\,\pi$;
\item[\emph{(v)}] \quad If $f$ maps $S$ onto $T$, then
$\varphi$ is an apartness isomorphism.
\end{itemize}
\end{theorem}

The next theorem deals with a more general situation.
\begin{theorem}  \label{ch-msa-cmrth3}
Let $S$ be a semigroup with apartness. Then
\begin{itemize}
\item[\emph{(i)}] \quad Let $\mu$ be a congruence, and $\kappa$ a co-congruence on $S$. Then, $\kappa$
  defines an apartness  on the factor set $S/\mu$ if and only if $\mu\cap\kappa =
  \emptyset$.
\item[\emph{(ii)}]  \quad The quotient mapping $\pi : S \rightarrow
S/\mu$, defined by $\pi (x)  = x\mu$, is an  se-epimorphism.
\end{itemize}
\end{theorem}

As a consequence of Theorem~\ref{ch-msa-basicfactor} and
Theorem~\ref{ch-msa-cmrth3} we have the following generalization of
Theorem~\ref{ch-msa-cmrth22}.
\begin{theorem} \label{ch-msa-basicfactor1}
  Let $f:S\to T$ be a mapping between sets with aparteness,  and let
  $\kappa$ be a co-equivalence on $S$ such that $\kappa\cap \mathrm{ker}\,f=\emptyset$.
Then:
\begin{itemize}
\item[\emph{(i)}] \quad  If $S$ is semigroup with apartness and $\kappa$ a
co-congruence, then $S/\mathrm{ker}\,f$ is a semigroup with
apartness, and $\pi:S\rightarrow S/\mathrm{ker}\,f $ an
se-epimorphism;
\item[\emph{(ii)}] \quad If, in addition, $T$ is a semigroup with apartness and $f$ an
se-homomorphism, then $\varphi:S/\mathrm{ker}\,f\rightarrow T$ is
also an se-homomorphism.
\end{itemize}
\end{theorem}

\subsubsection{Co-ordered semigroup with apartness}\label{ch-msa-324}

 The next results to be presented are consequences of the ones given
in Subsection~\ref{ch-msa-315} and Subsection~\ref{ch-msa-322}.

A tuple $(S,=,\#; \cdot; \alpha)$ is called
\begin{itemize}
\item[$\bullet$] \emph{co-quasiordered semigroup} if $\alpha$ is a
co-compatible co-quasiorder on $S$;
\item[$\bullet$] co-ordered semigroup if $\alpha$ is a
co-compatible co-order on $S$.
\end{itemize}

A constructive version for the Birkhoff's lemma for semigroups
follows.

\begin{lemma}\label{ch-msa-cb2}
Let $(S,=,\#; \cdot;\, \tau)$ be a co-quasiordered semigroup with
apartness. Then
 \begin{itemize}
\item[\emph{(i)}]\quad $\kappa_{\tau} = \tau\cup \tau^{-1}$ is a
  co-congruence on $S$;
\item[\emph{(ii)}] \quad  $(S/\kappa_{\tau}^c,=,\#)$ is a semigroup with apartness
where
\begin{eqnarray*}
x\kappa^c_\tau\,=y\kappa^c_\tau\  & \Leftrightarrow\ (x,y)\bowtie \kappa_\tau\\
x\kappa^c_\tau\,\#\,y\kappa^c_\tau\  & \Leftrightarrow\
(x,y)\in\kappa_\tau;
\end{eqnarray*}
\item[\emph{(iii)}] \quad $(S/\kappa^c_{\tau},=,\#;\, \Upsilon_\tau)$ is a
co-ordered semigroup with apartness where
$$(x\kappa^c_\tau, y\kappa^c_\tau)\in\Upsilon_\tau  \stackrel{\rm
def}{\Leftrightarrow} (x,y)\in\tau;$$
\item[\emph{(iv)}] \quad The quotient mapping $\pi_\tau: S\rightarrow
S/\kappa^c_\tau$, defined by $\pi_\tau(x)= x\kappa^c_\tau$, is an
isotone and reverse isotone se-epimorphism.
\end{itemize}
 \end{lemma}

\begin{theorem}\label{ch-msa-cp12}
Let $(S,=,\#)$ be a semigroup with apartness, and let $(T,=,\#;\,
\sigma)$ be a co-quasiordered semigroup with apartness. If
$f:S\rightarrow T$ is an se-mapping, then
 \begin{itemize}
\item[\emph{(i)}]\quad  $\mu_f$ is a co-compatible co-quasiorder on $S$;
\item[\emph{(ii)}] \quad $\kappa_f=\mu_f\cup\mu_f^{-1}$ is a co-congruence on on $S$ such that
$\kappa_f\subseteq \mathrm{coker}\,f$;
\item[\emph{(iii)}]\quad  If $\sigma$ is a co-order on $T$, then
$\kappa_f=\mathrm{coker}\,f$;
\item[\emph{(iv)}] \quad $(S/\kappa^c_f, =, \#; \cdot; \, \Upsilon_f)$ where
$$(x\kappa^c_f, y\kappa^c_f)\in\Upsilon_f  \stackrel{\rm
def}{\Leftrightarrow} (x,y)\in\mu_f$$ is a co-ordered semigroup with
apartness;
\item[\emph{(v)}] \quad The mapping $\psi: S/\kappa^c_f \rightarrow S/\kappa^c_\sigma$ between
co-ordered semigroups with apartness $(S/\kappa^c_f, =, \#; \cdot;
\, \Upsilon_f)$ and $(S/\kappa^c_\sigma , =, \#; \cdot\,
\Upsilon_\sigma)$, defined by $\psi(x\kappa_f)=
(f(x))\kappa_\sigma$, is an isotone and reverse isotone a-injective
se-homomorphism such that $\psi\pi_f=\pi_\sigma f$.
\end{itemize}
 \end{theorem}

\begin{theorem}
Let $(S,=,\#; \cdot\,\, \tau)$ be a co-quasiordered semigroup with
apartness, and let $(T,=,\#; \cdot; \, \sigma)$ be a co-quasiordered
semigroup with tight apartness. If $f:S\rightarrow T$ is an
se-homomorphism such that $\mu_f\subseteq\tau$, then the mapping
$\varphi: S/\kappa_\tau^c \rightarrow T$, defined by
$\varphi(x\kappa_\tau^c)=f(x)$, is an se-homomorphism such that
$\varphi\pi_\tau = f$.
\end{theorem}

\subsection{Applications and possible
applications}\label{ch-msa-33}

There is no  doubt about the deep connections between constructive
mathematics and computer science. Moreover, ``if programming is
understood  not as the writing of instructions for this or that
computing machine but as the design of methods of computation that
is the computer's duty to execute, then it no longer seems possible
to distinguish the discipline of programming from constructive
mathematics'', Martin-L\"{o}f  (1982)

Let us give some examples of applications of  ideas presented in the
previous section. We will start with constructive analysis. The
proof of one of the directions of the constructive version of the
Spectral Mapping Theorem is based on some elementary constructive
semigroups with inequality techniques, Bridges and Havea (2001).
 It is also worth mentioning the applications of
commutative basic algebraic structures with tight apartness within
the automated reasoning area, Calder\'{o}n (2017). For
possible applications within computational linguistic see Moshier (1995).
 Some topics from mathematical economics can be
approached constructively too (using some order theory for sets with
apartness), Baroni and Bridges (2008).

The theory of semigroups with apartness is, of course, in its
infancy, but, as we have already pointed out, it promises the
prospect of applications in other (constructive) mathematics
disciplines, certain areas of computer science, social sciences,
economics. Contrary to the classical case, the applications of
constructive semigroups with apartness, due to their novelty,
constitute an unexplored area.

To conclude, although one of the main motivators for initiating and
developing the theory of semigroups with apartness comes from the
computer science area, in order to have profound applications, a
certain amount of the theory, which can be applied, is necessary
first. Among priorities, besides growing the general theory, are
further developments of: constructive relational structures -
(co)quotient structures, constructive order theory, theory of
ordered semigroups with apartness, etc.

The study of basic constructive algebraic structures with apartness
as well as constructive algebra as a whole can impact the
development of other areas of constructive mathematics. On the other
hand, it can make both proof engineering and programming more
flexible.

\subsection{Notes}\label{ch-msa-34}

Constructive mathematics is not a unique notion. Various forms of
constructivism have been developed over time. The principal trends
include the following varieties: \textbf{INT} - Brouwer's
intuitionistic mathematics, \textbf{RUSS} - the constructive
recursive mathematics of the Russian school of Markov, \textbf{BISH}
- Bishop's constructive mathematics. Every form  has intuitionistic
logic at its core. Different schools have different additional
principles or axioms given by the particular approach to
constructivism. For example, the notion of  an \emph{algorithm} or a
\emph{finite routine} is taken as primitive in \textbf{INT} and
\textbf{BISH}, while \textbf{RUSS} operates with a fixed programming
language and an algorithm is a sequence of symbols in that language.
The Bishop-style of constructive mathematics enables one to
interpret the results both in classical mathematics, \textbf{CLASS},
and other varieties of constructivism. We regard classical
mathematics as Bishop-style mathematics plus the law of excluded
middle, \textbf{LEM}.

We have already emphasized that the Errett Bishop - style
constructive mathematics, \textbf{BISH},  forms the framework for
our work. \textbf{BISH} originated in 1967 with the publication of
the book \emph{Foundations of Constructive Mathematics},
 and with its second, much revised edition (Bishop and Bridges, 1985).
 There has been a steady stream of publications
contributing to Bishop's programme since 1967.  A ten-year long
systematic research of computable topology, using apartness as the
fundamental notion,  resulted in the first book, Bridges and
V\^{\i}\c{t}\={a} \cite{ch-msa-dbv2} (2011) on topology within
\textbf{BISH} framework. \texttt{Modern algebra}, as is noticed,
``contrary to Bishop's expectations, also proved amenable to a
natural, thoroughgoing, constructive treatment''.

Constructive algebra is a relatively old discipline developed among
others by L. Kronecker, van der Waerden, A. Heyting. For more
information on the history see Mines, Richman and Ruitenburg's book (1988).
 Troelstra and van Dalen's book (1988).
 One of the main topics in constructive algebra is constructive
algebraic structures with the relation of (tight) apartness $\#$,
the second most important relation in constructive mathematics. The
principal novelty in treating basic algebraic structures
constructively is that (tight) apartness becomes a fundamental
notion. (Consider the reals: we cannot assert that $x^{-1}$ exists
unless we know that $x$ is apart from zero, i.e. $|x| >0$ -
constructively that is not the same thing as $x\# 0$. Furthermore,
in fields $x^{-1}$ exists only if $x$ is apart from 0,
(Beeson, 1985)).

In some books and papers,  such as Troelstra and van Dalen's book  (1988),
 the term ``preapartness'' is used for an apartness
relation, while ``apartness'' means tight apartness. The tight
apartness on the real numbers was introduced by  Brouwer in the
early 1920s. Brouwer introduced the notion of apartness as a
positive intuitionistic basic concept. A formal treatment of
apartness relations began with A. Heyting's formalization of
elementary intuitionistic geometry  in  (Heyting, 1927). The
intuitionistic axiomatization of apartness is given in
 (Heyting, 1956).
 The study of algebraic structures in the presence of \emph{tight}
apartness  was started by  Heyting (1925). Heyting gave
the theory a firm base in 1941 (Heyting, 1941).

During the implementation of the FTA Project, Geuvers  {\it et al} (2002a)
 the notion of commutative constructive
semigroups with tight apartness appeared. We put noncommutative
constructive semigroups with ``ordinary'' apartness in the center of
our study, proving first, of course, that such  semigroups do exist.
Starting our work on constructive semigroups with apartness, as
 pointed out above,  we  faced   an algebraically completely new
area. What we had in ``hand'' at that moment were the experience and
knowledge coming from classical semigroup theory, other constructive
mathematics disciplines, and  computer science.

The notion of co-quasiorder first appeared in Romano (1996).
 However, let us mention that the results
reported from  (Romano, 1996): Theorem 0.4, Lemma 0.4.1, Lemma
0.4.2, Theorem 0.5 and Corollary 0.5.1 (pages 10-11 in (Romano, 2002))
are not correct. Indeed, the filled product
mentioned in Romano (1996, 2002) is
not associative in general. The notion of co-equivalence, i.e. a
symmetric co-quasiorder, first appeared in Bo\v zi\'c and Romano
(1987).

The  Quotient Structure Problem, \textbf{QSP},  is one of the very
first problems which has to be considered for any structure with
apartness. The solutions of the \textbf{QSP} problem for sets and
semigroups with apartness was given in Crvenkovi\'c  {\it et al} (2013).
 Those results are improved in Mitrovi\'c   {\it et al} (2019).

More background on constructive mathematics can be found in the
following books: Beeson (1985),
 Bishop (1967),
Bridges and  V\^{\i}\c{t}\={a} (2011), Troelstra and van
Dalen (1988). The standard reference for constructive
algebra are Mines,  Richman and  Ruitenburg (1988),
Ruitenburg (1982).

\section{\textbf{CLASS} and \textbf{BISH} -
a comparative analysis}\label{ch-msa-4}

From the classical mathematics (\textbf{CLASS}) point of view,
mathematics consists of a preexisting mathematical truth.  From a
constructive viewpoint,  the judgement $\varphi$ \emph{is true}
means that \emph{there is a proof of} $\varphi$.
 ``What constitutes  a proof is a social construct, an agreement
among people as to what is a valid argument. The rules of logic
codify a set of principles of reasoning that may be used in a valid
proof. Constructive (intuitionistic) logic codifies the principles
of mathematical reasoning as it is actually practiced,'' Harper (2013).
 In constructive mathematics, the \emph{status of
an existence statement} is much stronger than in \textbf{CLASS}. The
classical interpretation is that an object exists if its
non-existence is contradictory. In constructive mathematics when the
existence of an object is proved, the proof also demonstrates how to
find it. Thus, following further Harper (2013), the
\emph{constructive logic} can be described as the logic of people
matter, as distinct from the classical logic, which may be described
as the logic of the mind of God. One of the main features of
constructive mathematics is that the concepts that are equivalent in
the presence of \textbf{LEM}, need not be equivalent any more. For
example, we distinguish nonempty and inhabited sets, several types
of inequalities, two complements of a given set, etc.

 Contrary to the classical case, a set exists only when it is
defined. To define a set  $S$, we have to give a property that
enables us to construct members of $S$, and to describe the equality
$=$ between elements of $S$ -- which is a matter of convention,
except that it must be an  equivalence... There is another problem
to face   when we consider families of sets that are closed under a
suitable operation of complementation. Following Bishop and Bridges (1985)
 ``we do not wish to define complementation in the
terms of negation; but on the other hand, this seems to be the only
method available. The way out of this awkward position is to have a
very flexible notion based on the concept of \emph{a set with
apartness}.''

In \textbf{CLASS},   equivalence is the natural generalization of
equality.  A theory with equivalence involves
 equivalence and functions, and relations respecting this equivalence.
In constructive mathematics the same   works without difficulty,
Ruitenburg  (1991).

 Many sets  come with a binary relation called inequality satisfying certain properties, and
denoted by $\neq$, $\#$ or $\not\backsimeq$. In general, more
computational information is required to distinguish elements of a
set $S$, than to show that elements are equal. Comparing with
\textbf{CLASS}, the situation for inequality
 is more complicated. There are different types of
inequalities (denial inequality, diversity, apartness, tight
apartness - to mention a few), some of them completely independent,
which only in \textbf{CLASS} are equal to one standard inequality.
So, in \textbf{CLASS} the study of the equivalence relation
suffices, but in constructive mathematics, an inequality becomes a
``basic notion in intuitionistic axiomatics''.  Apartness, as a
positive version of inequality, in the words of Jacobs (1995), ``is yet another fundamental notion developed in
intuitionism which shows up in computer science.''

The statement that every equivalence relation is the negation of
some apartness relation is equivalent to the excluded middle. The
statement that the negation of an equivalence relation is always an
apartness relation is equivalent to the nonconstructive de Morgan
law.

For a tight apartness, the two complements are constructive
counterparts of the classical complement. In general, we have ${\sim
Y} \subseteq {\neg Y}.$ However, even for a tight apartness, the
converse inclusion entails the Markov principle, \textbf{MP}. This
result illustrates a main feature of constructive mathematics:
classically equivalent notions could be no longer equivalent
constructively. For which  type of subset  of a set with apartness
do we  have equality between its two complements?
 It turns out that the answer initiated a development of \emph{order theory} for sets and
 semigroups with apartness we develop. Constructive mathematics brings
 to the light some notions which are invisible to the classical eye
 (here, the three notions of detachability).

In the constructive order theory, the notion of co-transitivity,
that is the property that for every pair of related elements, any
other element is related to one of the original elements in the same
order as the original pair is a constructive counterpart to
classical transitivity, (Crvenkovi\'c  {\it et al}, 2013).

A relation defined on a set with apartness $S$ is
\begin{itemize}
\item[$\bullet$] \emph{weak co-quasiorder} if it is irreflexive and
cotransitive,
\item[$\bullet$] \emph{co-quasiorder } if it is strongly irreflexive and
cotransitive.
\end{itemize}
Even if the two classically   (but not constructively) equivalent
variants of a co-quasiorder are constructive counterparts of a
quasiorder in the case of (a tight) apartness, the stronger variant,
co-quasiorder, is, of course, the most appropriate for a
constructive development of the theory of semigroups with apartness
we develop. The weaker variant, that is, weak co-quasiorder, could
be relevant in analysis. ``One  might expect that the splitting of
notions leads to an enormous proliferation of results in the various
parts of constructive mathematics when compared with their classical
counterparts. In particular, usually only very few constructive
versions of a classical notion are worth developing since other
variants do not lead to a mathematically satisfactory theory,''
Troelstra and van Dalen  (1988).

Following Bishop, every classical theorem presents the challenge:
find a constructive version with a constructive proof. Within
\textbf{CLASS}, the semigroups can be viewed, historically, as an
algebraic abstraction of the properties of the composition of
 transformations on a set. Cayley's theorem for semigroups (which
can be seen as an extension of the celebrated Cayley's theorem on
groups) stated that every semigroup can be embedded in a semigroup
of all self-maps on a set. As a consequence of
Theorem~\ref{ch-msa-evid}, we can formulate \emph{the constructive
Cayley's theorem for semigroups with apartness}.

Following the standard literature on constructive mathematics, the
term ``constructive theorem'' refers to a theorem with a
constructive proof. A classical theorem that is proven in a
constructive manner is a constructive theorem. This constructive
version can be obtained by strengthening the conditions or weakening
the conclusion of the theorem. Although constructive theorems might
look like the corresponding classical versions, they often have more
complicated hypotheses and proofs. Theorems and their proofs given
in
\begin{itemize}
\item[$\bullet$] subsections 2.1.2, 2.1.3, 2.2.2 and 2.2.3 for
classical case,
\item[$\bullet$] subsections 3.1.4, 3.1.5, 3.2.3 and 3.2.4 for
constructive case,
\end{itemize}
are evidence for that.

There are, often, several constructively different versions of the
same classical theorem.  Some  classical theorems are neither
provable nor disprovable, that is, they are independent of
\rm{\textbf{BISH}}. For some classical theorems it is shown that
they are not provable constructively. More details about
nonconstructive principles and various classical theorems that are
not constructively valid can be found in Ishihara (2013).

\section{Concluding remarks}\label{ch-msa-5}

There are many decisions a mathematician must make when deciding to
replace classical logic with intuitionisitic logic. Let us mention
some of them. First of all - \textbf{choice of variant of
constructive mathematics.} Constructive mathematics is not a unique
notion. Our choice was the Errett Bishop - style constructive
mathematics, \textbf{BISH}. The cornerstones for \textbf{BISH}
include the notion of positive integers, sets and functions.
Contrary to the classical case, a set exists only when it is
defined.

Going through the literature there are several variants of what is
considered to be a set with apartness - depending on the relations
between equality and apartness defined on a set. Our choice - our
starting structure - is a \textbf{set with apartness} $(S, =, \#)$
where
\begin{itemize}
\item[$\bullet$] equality and apartness are  basic notions,
\item[$\bullet$] equality and apartness are independent of each
other,
\item[$\bullet$] apartness is not, in general, tight.
\end{itemize}
Such our choice - fully justified in Darp\"{o} and Mitrovi\'c
(arXiv:2103.07105) - was and is a novelty within constructive circles.

When working constructively one has to choose definitions with some
care. Heyting pointed out that it is very important to have
experience of where the intuitionistic pitfalls lie. Classically
equivalent notions often split into a certain number of inequivalent
constructive ones: nonempty and inhabited set, several inequivalent
constructive definition of inequality; two complements of a given
subset of a set - just to mention a few. Thus,  it is very important
to be aware of the phenomena of splitting notions.

Is it easier to work within constructive mathematics than within
classical? One does not develop constructive mathematics  for
simplification. There is something else, something so important
worth  that simplicity can  sometimes be sacrificed. So WHY to work
constructively? Is it usefulness in question? Or, something else?
Instead of giving answers, let us cite Heyting, {\it Intuitionism - An
Introduction} (1956)  again.

``It seems quite reasonable to judge a mathematical system by its
usefulness [...]  in my eyes its chances of being useful for
\textbf{philosophy, history and the social sciences} are better. In
fact, mathematics, from the intuitionistic point of view, is a study
of certain functions of the human mind, and as such it is akin to
these sciences. But is usefulness really the only measure of value?
You know how philosophers struggle with the problem of defining the
concept of value in art; yet every educated person feels this value.
The case is analogous for the value of intuitionistic mathematics.''

\section*{Acknowledgements} The authors are grateful
to anonymous referees for careful reading of the manuscript and
helpful comments. M. M.  was financially supported by the Ministry
of Science, Technological Development and Innovation of the Republic
of Serbia (Contract No. 451-03-47/2023-01/ 200109). M. N. H. is
supported by the ICMPA-UNESCO Chair, which  is in partnership with Daniel
Iagolnitzer Foundation (DIF), France, and the Association pour la
Promotion Scientifique de l'Afrique (APSA), supporting the
development of mathematical physics in Africa. P. C. also thanks
the Portuguese Foundation for Science and Technology (FCT) under the projects
UIDB/00013/2020 and UIDP/00013/2020.

\section{Appendix}

\subsection{The \textbf{CLASS} case}

\begin{proof} \ - \textbf{Lemma~\ref{msart-11o}}

Straightforward, based on the well-known facts that each element of
$S$ belongs to a unique equivalence class,  and each equivalence
class is nonempty.
\end{proof}

\begin{proof}\  - \textbf{Theorem~\ref{msart-11}}

(i).  Let $x(ker\,f) = y(ker\,f)$. Then $(x,y)\in ker\,f$, thus
$f(x)=f(y)$ and $\varphi$ is well-defined, i.e. the definition of
$\varphi(x(kef\,f))$ is independent of the choice of element form
$x(ker\,f)$.

Let $x\in S$. By the assumption
$$(\varphi\circ\pi_{nat}) (x)= \varphi(\pi_{nat} (x)) = \varphi(x(ker\,f))= f(x).$$
If $\varphi(x(ker\,f)) =\varphi(y(ker\,f))$ then $f(x)=f(y)$, which,
further, means that $(x,y)\in ker\,f$. Thus $x(ker\,f) = y(ker\,f)$,
and $\varphi$ is injective.

If there is another function $\varphi' : S/ker\,f \rightarrow T$
such that $f=\varphi'\circ\pi_{nat}$, then for any $x(ker\,f)\in
S/ker\,f$ we have
$$\varphi'(x(ker\,f)= \varphi'(\pi_{nat}(x))= \varphi'\circ\pi_{nat}(x)= f(x).$$
Thus $\varphi=\varphi'$.

(ii). If $f$ is an  onto mapping then for any $y\in T$ there is an
$x\in S$ such that $y=f(x)$. But then $y=\varphi(x(ker\,f))$ and so
$\varphi$ is surjective. We have that $\varphi$ is bijection.

(iii). Let $\varphi$ be a surjective mapping. By
Lemma~\ref{msart-11o}, $\pi_{nat}$ is surjective too. Thus
$\varphi\circ\pi$ is also surjective, i.e. $f=\varphi\circ\pi_{nat}$
is surjective.
\end{proof}

\begin{proof} \  - \textbf{Theorem~\ref{msart-12}}

Let $x\varepsilon = y\varepsilon$, $x,y\in S$. Then $(x,y)\in
\varepsilon\subseteq ker\,f$. Thus $f(x)=f(y)$ and $\varphi$ is well
defined. It is a routine to verify that
$f=\varphi\circ\pi_{\varepsilon}$. The rest follows by the
Theorem~\ref{msart-11}.
\end{proof}

\begin{proof} \ - \textbf{Lemma~\ref{pl1}}

(i). \ Let  $f$  be an isotone mapping. If $(x,y)\in \rho$, then, by
the assumption, $(f (x), f(y))\in \sigma$, which, by the definition
of $\eta_f$ given above, means that $(x,y)\in \eta_f$.

Conversely, from the assumption and the definition of $\eta_f$ we
have
$$(x,y)\in\rho \Rightarrow (x,y)\in \eta_f
\Leftrightarrow (f (x), f(y))\in \sigma.$$

(ii). \ Similar to (i).

(iii). Consequence of (i) and (ii).
\end{proof}

\begin{proof} \ - \textbf{Lemma~\ref{rel-structure}}

Straightforward.
\end{proof}

\begin{proof} \ -  \textbf{Lemma~\ref{msart-bl}}

(i).\ Follows immediately by Lemma~\ref{msart-14}

(ii). \ By  Lemma~\ref{rel-structure}, $(S/\varepsilon_{\rho};
\hat{\Theta}_\rho)$ ia a quasiordered set. Let  $(x\varepsilon_\rho,
y\varepsilon_\rho), (y\varepsilon_\rho, x\varepsilon_\rho)
\in\hat{\Theta}_\rho$, which, by the definition, gives $(x,y),
(y,x)\in \rho$, i.e. $(x,y)\in \rho\cap\rho^{-1}=\varepsilon_\rho$.
Thus $x\varepsilon_\rho = y\varepsilon_\rho$. So,
$\hat{\Theta}_\rho$ is an order.

(iii). Let $x,y\in  S$. Then
\begin{align*}
(x,y)\in\rho &\Leftrightarrow ( x\varepsilon_{\rho},
y\varepsilon_{\rho})\in \Theta_{\rho}\\
&\Leftrightarrow (\pi_{nat}(x), \pi_{nat}(y))\in \Theta_{\rho}.
\end{align*}
Thus $\pi$ is isotone and reverse isotone, and by the
Lemma~\ref{msart-11o}, $\pi$ is surjective.
\end{proof}

\begin{proof} \ -  \textbf{Lemma~\ref{msart-15}}

(i)$\Rightarrow$(ii). \   From $(x,x)\in\varepsilon_{\rho}$ we have
$x\in x\varepsilon_{\rho}$. Also $y\in y\varepsilon_{\rho}$. Thus,
by (i) we have $(x,y)\in\rho$.

(ii) $\Rightarrow$ (i). \ Follows immediately from $x\in
x\varepsilon_{\rho}$ and $y\in y\varepsilon_{\rho}$.

(i) $\Rightarrow$ (iii). Let $x_1\in x\varepsilon_{\rho}$ and
$y_1\in y\varepsilon_{\rho}$,i.e. we have $(x_1,y_1)\in\rho$. Let
$a\in x\varepsilon_{\rho}$ and $b\in y\varepsilon_{\rho}$. Since
$(a,x),(x_1,x)\in \varepsilon_{\rho}$ we have $(a,x_1)\in
\varepsilon_{\rho}$. Since $(y,y_1),(y,b)\in\varepsilon_{\rho}$, we
have $(y_1,b)\in\varepsilon_{\rho}$. Then $(y_1,b)\in\rho$. Now,
from $(a,x_1),(x_1,y_1), (y_1,b)\in\rho$ we have $(a,b)\in\rho$.

(iii) $\Rightarrow$ (i). Follows immediately.
\end{proof}

\begin{proof} \  -  \textbf{Theorem~\ref{d-msart-11}}

(i).\ Reflexivity of $\sigma$ gives $(f(x), f(x))\in\sigma$, $x\in
S$, and, by the definition of $\eta_f$, we have $(x,x)\in\eta_f$.

Transitivity of $\eta_f$ follows by the definition of $\eta_f$ and
transitivity of $\sigma$, i.e. we have
\begin{align*}
(x,y),(y,z)\in\eta_f&\Leftrightarrow
(f(x),f(y)),(f(y),f(z))\in\sigma\\
&\Rightarrow (f(x),f(z))\in\sigma\\
&\Leftrightarrow (x,z)\in\eta_f.
\end{align*}

(ii).\ By (i) and Lemma~\ref{msart-bl}(i), $\varepsilon_f
=\eta_f\cap \eta_f^{-1}$ is an equivalence on $S$.

Using the reflexivity of $\sigma$ we have

\begin{align*}
(x,y)\in ker\,f &\Leftrightarrow f(x)= f(y)\\
&\Rightarrow (f(x),f(y))\in\sigma \\
&\Leftrightarrow (x,y)\in\eta_f.
\end{align*}
i.e. $ker\,f \subseteq \eta_f$. In a similar manner we can prove
$ker\, \subseteq \eta_f^{-1}$. Thus  $ker\, f\subseteq
\varepsilon_f$. Conversely,let
\begin{align*}
(x,y)\in \varepsilon_f =\eta_f\cap \eta_f^{-1} &\Leftrightarrow
(x,y), (y,x)\in \eta_f \\
&\Leftrightarrow (f(x),f(y)), (f(y),f(x))\in\sigma \\
& \Rightarrow f(x)= f(y)\\
&\Leftrightarrow (x,y)\in ker\,f,
\end{align*}
so $\varepsilon_f \subseteq ker\,f$. Thus,  $ker\, f
=\varepsilon_f$.

(iii).\  By Theorem~\ref{msart-11}(i), the mapping $\varphi :
S/\varepsilon_f \rightarrow T$ defined by $\varphi(x\varepsilon_f))
= f(x)$ is the unique injective mapping such that
$f=\varphi\circ\pi_{nat}$.

By Lemma~\ref{msart-bl},  $(S/\varepsilon_f, \Theta_f)$ is an
ordered set, and
\begin{align*}
(x\varepsilon_f,y\varepsilon_f)\in  =\Theta_f &\Leftrightarrow
(x,y)\in \eta_f\\
&\Leftrightarrow (f(x),f(y))\in\sigma.
\end{align*}
Therefore, $\varphi$ is isotone and reverse isotone  mapping such
that $f=\varphi\circ\pi$.

(iv). \ Follows by (ii) and Theorem~\ref{msart-11}(iii).
\end{proof}

\begin{proof} \  -  \textbf{Theorem~\ref{d-msart-112}}

(i).\ From $\rho\subseteq\eta_f$ we have
$\rho^{-1}\subseteq\eta_f^{-1}$. Then
$$\varepsilon_{\rho}= \rho\cap\rho^{-1}\subseteq
\eta_f\cap\eta_f^{-1}=\varepsilon_f.$$ Now, by the
Theorem~\ref{d-msart-11}, we have
$\varepsilon_{\rho}\subseteq\varepsilon_f= ker\, f$.

(ii) By the Theorem~\ref{msart-12} there is the unique injective
mapping $\varphi : S/\varepsilon_{\rho} \rightarrow T$, defined by
$\varphi(x\varepsilon_{\rho}) = f(x)$, such that
$f=\varphi\circ\pi_{\rho}$. By Lemma~\ref{msart-bl},
$(S/\varepsilon_{\rho},\hat{\Theta}_\rho)$ ia an ordered set.

Let $(x\varepsilon_{\rho}, y\varepsilon_{\rho})\in
S/\varepsilon_{\rho}$. Then
\begin{align*}
( x\varepsilon_{\rho}, y\varepsilon_{\rho})\in \hat{\Theta}_\rho
&\Leftrightarrow (x,y)\in \rho\subseteq \eta_f\\
&\Rightarrow ( f(x), f(y))\in \sigma.
\end{align*}
So, $\varphi$ is an isotone mapping.
\end{proof}

\begin{proof} \  - \textbf{Lemma~\ref{ch-msart-3-l1}}

\ Let $\rho$ be a compatible quasioder on $S$. Let  $(x,y)\in \rho$
and $z\in S$. By the reflexivity of $\rho$, we have $(z,z)\in S$.
Now, by the compatibility, it follows that $(zx,zy)\in S$, and
$\rho$ is  left compatible. Right compatibility can be proved in a
similar manner.

Conversely, Let $\rho$ be a left and right compatible quasiorder on
$S$. Let $(x,y), (s,t)\in \rho$. By the left and right compatibility
it follows that $(xs,ys), (ys,yt)\in\rho$, which, by the
transitivity of $\rho$, gives us $(xs,yt)\in\rho$.
\end{proof}

\begin{proof} \  -  \textbf{Theorem~\ref{sth2-n}}

(i). \ All the set-theoretic results are as in
Theorem~\ref{msart-11}. It only remains to prove that $ker\,f$ is
compatible with multiplication as well as that $\varphi$ is a
homomorphism.

Let $(x,y), (s,t)\in ker\,f$. Then $f(x)=f(y)$ and $f(s)=f(t)$. It
follows that
$$f(xs)=f(x)f(s)=f(y)f(t)=f(yt),$$
and hence $(xs,yt)\in ker\,f$. Therefore $ker\,f$ is a congruence.

Let $x(ker\,f), y(ker\,f)\in  S/ker\,f$. Then
\begin{align*}
\varphi(x(ker\,f) y(ker\,f)) &= \varphi((xy)(ker\,f))\\
&=f(xy)=f(x)f(y)\\
&=\varphi(x(ker\,f)) \varphi(y(ker\,f)).
\end{align*}
So, $\varphi$ is a homomorphism.

(ii). \ Follows immediately by the assumptions and
Theorem~\ref{msart-12}.
\end{proof}

\begin{proof} \  -  \textbf{Lemma~\ref{msart-3-bl}}

(i).\ Let $(x,y),(s,t)\in\rho^{-1}$, or, equivalently
$(y,x),(t,s)\in \rho$. Then, as $\rho$ is compatible,
$(yt,xs)\in\rho$, which further means that $(xs,yt)\in\rho^{-1}$.
Thus $\rho^{-1}$ is a compatible quasiorder, and the equivalence
$\varepsilon_{\rho}= \rho\cap\rho^{-1}$ is compatible too, i.e. it
is a congruence on $S$.

(ii).\ Compatibility of $\hat{\Theta}_\rho$ follows by its
definition and compatibility of $\rho$.

(iii). \ By  Theorem~\ref{sth21} $\pi$ is an epimorphism. By
Lemma~\ref{msart-bl}(iii), it is isotone and reverse isotone.
\end{proof}

\begin{proof} \  - \textbf{Theorem~\ref{d-msart-3-11}}

Follows by  Theorem~\ref{d-msart-11} and Theorem~\ref{sth21}.
\end{proof}

\begin{proof} \  -  \textbf{Theorem~\ref{d-msart-3-111}}

Follows by  Theorem~\ref{d-msart-112} and
Theorem~\ref{d-msart-3-11}.
\end{proof}

\subsection{Within BISH}

\begin{proof} \  -  \textbf{Theorem~\ref{ch-msa-thm1.2}}

See (Mines,  Richman and  Ruitenburg, 1988,  Theorem I 2.2).
\end{proof}

\begin{proof} \  -  \textbf{Theorem~\ref{senegcom}}

(i). \ Let $Y$ be an sd-subset of $S$. Then, applying the definition
and  logical axiom we have
\begin{eqnarray*}
\forall _{ x\in S} \, (x \in Y \vee x \in {\sim} Y)\
&\Leftrightarrow \ &
\forall _{ x\in S} \,(x \in Y \, \vee \,  \forall _{ y\in Y} (x\# y))\\
 &\Rightarrow & \forall _{ x\in S} \,\forall _{ y\in Y} \, (x \in Y \vee x \# y).
\end{eqnarray*}

In order to prove the second part of this statement, we consider the
real number set $\mathbb{R}$ with the usual (tight) apartness and
the subset $Y=\widetilde{0}.$ Then, for each real number $x$ and for
each $y\in Y$ it follows, from the co-transitivity of $\#$, either
$y\# x$ or $x\# 0 $, that is, either $x\in Y$ or $x\# y$.
Consequently, $Y$ is a qd-subset of $\mathbb{R}.$ On the other hand,
if $Y$ is an sd-subset of $\mathbb{R},$ then for each $x\in
\mathbb{R}$, either $x\in Y$ or $x\in \sim Y.$ In the former case,
$x\# 0$ and in the latter $x=0$, hence \textbf{LPO} holds.

(ii). \ Let $Y$ be a qd-subset, and let $a\in \neg Y$. By assumption
we have
$$\forall _{ x\in S} \,\forall _{ y\in Y} \, (x
\in Y \vee x \# y),$$ so substituting $a$ for $x$, we get
$\forall_{y\in Y} \, (a\in Y \vee a\# y)$, and since, by assumption,
$\neg(a\in Y)$, it follows that $a\# y$ for all $y\in Y$. Hence
$a\in \sim Y$. See also \cite{ch-msa-mmscr}.

(iii). \ Let $S$ be the real number set $\mathbb{R}$ with the usual
 apartness $\#$. As in the proof of (i), consider the qd-subset $\widetilde{0}$
of $\mathbb{R}$. If $\widetilde{0}$ is a d-subset of $\mathbb{R}$,
then $x\in \widetilde{0}$ or $\neg(x\in \widetilde{0})$, for all
real numbers $x$. In the latter case $\neg(x\# 0)$,  which is
equivalent to $x=0.$ Thus we obtain the property $\forall_{ x\in
\mathbb{R}}\; (x\# 0\vee x=0)$ which, in turn, is equivalent to
\textbf{LPO}.

(iv).  \ Consider a real number $a$ with $\neg (a=0)$ and let $S$ be
the set $\{ 0,a\}$ endowed with the usual apartness of $\mathbb{R}$.
For $Y=\{ 0\}$, since $0\in Y$ and $a\in \neg Y,$ it follows that
$Y$ is  a d-subset of $S$. On the other hand, if $Y$ is a qd-subset
of $S$, then $a\# 0$. It follows that for any real number with $\neg
(a=0)$, $a\# 0$ which entails the Markov Principle, \textbf{MP}.

(v). \ The first part follows immediately from (i), (ii) and the
definition of d-subsets. The converse follows  from (i) and (iv).

(vi). \ Consider again $\mathbb{R}$ with the usual apartness and
define $Y=\{ 0\} .$ If $Y$ is a qd-subset of $\mathbb{R}$, then for
all $x\in \mathbb{R}$ we have  $x=0$ or $x\# 0$, hence \textbf{LPO}
holds.
\end{proof}

\begin{proof}   \  -  \textbf{Lemma~\ref{ch-msa-*properties}}

 If $\alpha\subseteq\gamma$ and $\beta\subseteq\delta$ then
\begin{eqnarray*}
(x,y)\in\alpha\ast\beta \  &\Leftrightarrow & \ \forall _{z\in S} \, ((x,z)\in\alpha \vee (z,y)\in\beta)\\
&\Leftrightarrow  & \ \forall _{z\in S} \, ((x,z)\in\gamma \vee (z,y)\in\delta)\\
&\Leftrightarrow & \ (x,y)\in\gamma\ast\delta.
\end{eqnarray*}
\end{proof}

\begin{proof} \  -  \textbf{Lemma~\ref{ch-msa-pl2}}

(i). \ If $\alpha$ is strongly irreflexive then
\begin{eqnarray*}
(x,y)\in\alpha\ast\beta \  &\Leftrightarrow & \ \forall _{z\in S} \, ((x,z)\in\alpha \vee (z,y)\in\beta)\\
&\Rightarrow  & \ (x,x)\in\alpha \vee (x,y)\in\beta\\
&\Rightarrow & \ (x,y)\in\beta.
\end{eqnarray*}
The case when $\beta$ is strongly irreflexive is analogous.

(ii). By (iv),  we have $\alpha \ast\alpha \subseteq \alpha\subseteq
\#$. Thus, $\alpha \ast\alpha$ is strongly irreflexive.
\end{proof}

\begin{proof} \  -  \textbf{Lemma~\ref{ch-msa-reflexiveconsistent}}

(i). \ Let $\alpha$ be a strongly irreflexive  relation on $S$. For
each $a\in S$, it can be easily proved that $(a,a)\# (x,y)$ for all
$(x,y)\in \alpha$.

Let $\sim\!\,\alpha$ be reflexive, that is $(x,x)\in\sim\!\,\alpha$,
for any $x\in S$. On the other hand, the definition of the
a-complement implies $(x,y)\# (x,x)$ for any $(x,y)\in \alpha$. So,
$x\# x$ or $x\# y$. Thus, $x\# y$, that is, $\alpha$ is strongly
irreflexive.

(ii). \ If $\alpha$ is symmetric, then
\begin{eqnarray*}
(x,y)\in\sim\!\,\alpha \ &\Leftrightarrow & \ \forall _{(a,b)\in \alpha} \, ((x,y)\#(a,b))\\
&\Rightarrow  & \ \forall _{(b,a)\in \alpha} \, ((x,y)\#(b,a))\\
 &\Rightarrow & \ \forall _{(b,a)\in \alpha} \, (x\# b \ \vee \ y\# a)\\
 &\Rightarrow & \ \forall _{(a,b)\in \alpha} \, ((y,x)\# (a,b))\\
&\Leftrightarrow & \ (y,x)\in\sim\!\,\alpha.
\end{eqnarray*}

(iii). \ If $(x,y)\in\sim\alpha$ and $(y,z)\in\sim\alpha$, then, by
the definition of $\sim\alpha$, we have that $(x,y)\bowtie \alpha$
and $(y,z)\bowtie \alpha$. For an element $(a,b)\in\alpha$, by
co-transitivity of $\alpha$, we have $(a,x)\in\alpha$ or
$(x,y)\in\alpha$ or $(y,z)\in\alpha$ or $(z,b)\in\alpha$. Thus
$(a,x)\in\alpha$ or $(z,b)\in\alpha$, which implies that $a\# x$ or
$b\# z$, that is $(x,z)\#(a,b)$. So, $(x,z)\bowtie \alpha$ and
$(x,z)\in\ \sim\alpha$. Therefore, $\sim\alpha$ is transitive.

(iv). \ Let $(x,y), (y,x)\in\sim\alpha$. Assume $x\# y$. Then, by
the co-antisymmetry of $\alpha$, we have $(x,y)\in \alpha$ or
$(y,x)\in \alpha$ which is impossible. Thus, $\neg(x\# y$), i.e.
$x=y$ as apartness is tight.
\end{proof}

\begin{proof} \  -  \textbf{Proposition~\ref{ch-msa-senegcompl}}

(i). \ Let $(x,y)\in S\times S$. Then, for all $(a,b)\in\tau$,
\begin{eqnarray*}
    a\tau x \vee x\tau b  &\ \Rightarrow\ & a\tau x \vee x\tau y \vee y\tau b\\
  & \ \Rightarrow\ & a\# x \vee x\tau y\vee y\# b \\
  & \ \Rightarrow\ & (a,b)\# (x,y) \vee x\tau y,
\end{eqnarray*}
that is, $\tau$ is a qd-subset.

(ii). \ It follows from (i) and Theorem~\ref{senegcom}(ii).

(iii). \  This is a consequence of
Lemma~\ref{ch-msa-reflexiveconsistent}(i), (iii).
\end{proof}

\begin{proof} \  -  \textbf{Lemma~\ref{ch-msa-cun}}

From the strong irreflexivity of $\tau$  and $\sigma$, i.e. from
$\tau\subseteq \#$ and $\sigma\subseteq \#$, we have $\tau
\cup\sigma\subseteq \#$. Thus, $\tau \cup\sigma$ is strongly
irreflexive.

 From $\tau \subseteq \tau\cup\sigma$ and
$\sigma\subseteq\tau\cup\sigma$, by Lema~\ref{ch-msa-*properties} we
have
\begin{eqnarray*}
\tau \ast \tau\,&\subseteq &\, (\tau\cup\sigma)\ast (\tau\cup\sigma)\\
\sigma\ast\sigma \,&\subseteq &\, (\tau\cup\sigma)\ast
(\tau\cup\sigma).
\end{eqnarray*}
By the assumption and by the Proposition~\ref{ch-msa-senegcompl}(i),
we have
\begin{eqnarray*}
\tau \,&\subseteq &\, (\tau\cup\sigma)\ast (\tau\cup\sigma)\\
\sigma \,&\subseteq &\, (\tau\cup\sigma)\ast(\tau\cup\sigma),
\end{eqnarray*}
which gives us $\tau\cup\sigma \subseteq (\tau\cup\sigma)\ast
(\tau\cup\sigma)$. Thus $\tau\cup\sigma$ is  co-transitive.
\end{proof}

\begin{proof} \  -  \textbf{Lemma~\ref{ch-msa-l2b}}

Indeed, let $\alpha$ be a co-quasiorder and $\beta$ an equivalence
on $S$ such that  $\mathbf{\alpha}$ defines an apartness on
$\mathbf{S/\beta}$. Let $(x,a), (y,b) \in\beta$, i.e. $a\in x\beta$
and $b\in y\beta$, which, by the assumption, gives $a\beta = x\beta$
and $b\beta = y\beta$. If $(x,y)\in \alpha$, then, by (Ap6), $x\beta
\,\#\, y\beta$, which, by (Ap5'), gives $a\beta \,\#\, b\beta$. By
(Ap6) we have $(a,b)\in\alpha$. In a similar manner, starting from
$(a,b)\in\alpha$ we can conclude $(x,y)\in \alpha$.
\end{proof}

\begin{proof} \  - \textbf{Theorem~\ref{ch-msa-t-coe}}

(i). \ By the Proposition~\ref{ch-msa-senegcompl}(iv),  $\kappa^c$
is a quasiorder on $S$, and, by the
Lemma~\ref{ch-msa-reflexiveconsistent}(ii), it is symmetric. Thus,
$\kappa^c$ is an equivalence.

Let $(x,y)\in\kappa$ and $(y,z)\in\kappa^c$. Thus $(x,y)\in \kappa$
and $(y,z)\bowtie \kappa$. By the symmetry and co-transitivity of
$\kappa$ we have $(x,z)\in\kappa$ or $(y,z)\in\kappa$. Thus
$(x,z)\in\kappa$, and $\kappa\looparrowleft \kappa^c$.

(ii).\  The strong irreflexivity of $\#$ is implied by its
definition and by the strong irreflexivity of $\kappa$.

Let $a\kappa^c\,\# \,b\kappa^c$. Then $(a,b)\in\kappa$ implies that
$(b,a)\in\kappa$, that is $b\kappa^c\,\# \,a\kappa^c$.\label{ere}

Let $a\kappa^c\,\#\, b\kappa^c$ and $u\kappa^c\in S/\kappa^c$. Then
$(a,b)\in\kappa$, and, by the co-transitivity of $\kappa$, we have
$(a,u)\in\kappa$ or $(u,b)\in\kappa$. Finally we have that
$a\kappa^c\# u\kappa^c$ or $u\kappa^c\,\#\, b\kappa^c$, so the
relation $\#$ is co-transitive.

(iii). \ Let $\pi(x)\#\pi(y)$, i.e. $x\kappa^c\,\#\, y\kappa^c$,
which, by what we have just proved, means that $(x,y)\in\kappa$.
Then, by the strong irreflexivity  of $\kappa$, we have $x\# y$. So
$\pi$ is an se-mapping.

Let $a\kappa^c\in S/\sim\kappa$ and $x\in a\kappa^c$. Then
$(a,x)\in\sim\kappa$, i.e. $a\kappa^c=x\kappa^c$, which implies that
$a\kappa^c=x\kappa^c=\pi(x)$. Thus $\pi$ is an se-surjection.
\end{proof}

\begin{proof} \  - \textbf{Theorem~\ref{ch-msa-cmrth2}}

(i) \ The strong irreflexivity of $\mathrm{coker}\, f$ is easy to
prove: if $(x,y)\in\mathrm{coker}\,f$, then $f(x)\# f(y)$ and
therefore $x\# y$.

If $(x,y)\in\mathrm{coker}\,f$, then, by the symmetry of apartness
in $T$, $f(y)\# f(x)$; so $(y,x)\in\mathrm{coker}\,f$.

If $(x,y)\in\mathrm{coker}\,f$ and $z\in S$, i.e. $f(x)\# f(y)$ and
$f(z)\in T$, then either $f(x)\# f(z)$ or $f(z)\# f(y)$;
 that is, either $(x,z)\in\mathrm{coker}\,f$ or $(z,y)\in
\mathrm{coker}\,f$. Hence $\mathrm{coker}\,f$ is a co-equivalence on
$S$.

(ii) \ Let $(x,y)\in\mathrm{coker}\,f$ and
$(y,z)\in\mathrm{\ker\,f}$. Then $f(x)\# f(y)$ and $f(y)=f(z)$.
Hence $f(x)\# f(z)$, that is, $(x,z)\in\mathrm{coker}\,f$, and
$\mathrm{coker}\,f \looparrowleft \mathrm{ker}\,f$.

Now let $(x,y)\in\ker\,f$, so  $f(x)=f(y)$. If
$(u,v)\in\mathrm{coker}\,f$, then, by the co-transitivity of
$\mathrm{coker}\,f$, it follows that $(u,x)\in\mathrm{coker}\,f$ or
$(x,y)\in\mathrm{coker}\,f$ or $(y,v)\in\mathrm{coker}\,f$. Thus
either $(u,x)\in\mathrm{coker}\,f$ or $(y,v)\in \mathrm{coker}\,f$,
and, by the strong irreflexivity of $\mathrm{coker}\,f$, either $u\#
x$ or $y\# v$; whence we have $(x,y)\#(u,v)$. Thus $(x,y)\bowtie
\mathrm{coker}\,f$, or, equivalently $(x,y)\in\
(\mathrm{coker}\,f)^c$.

(iii) \ This follows from the definition of $\#$ in $S/\ker\,f$ and
(i).

(iv)\ Let us first prove that $\varphi$ is well defined. Let
$x(\ker\,f),y(\ker\,f)\in S/\ker\,f$ be such that
$x(\ker\,f)=y(\ker\,f\mathbf{)}$, that is, $(x,y)\in\ker\,f$. Then
we have $f(x)=f(y)$, which, by the definition of $\varphi$, means
that $\varphi (x(\ker\,f))=\varphi(y(\ker\,f))$.

Now let $\varphi(x(\ker \,f))=\varphi(y(\ker\,f))$; then
$f(x)=f(y)$. Hence $(x,y)\in\ker\,f$, which implies that
$x(\ker\,f)=y(\ker\,f)$. Thus $\varphi$ is an injection.

Next let $\varphi(x(\ker\,f))\#\varphi(y(\ker\,f))$; then $f(x)\#
f(y)$. Hence $(x,y)\in\mathrm{coker}\,f$, which, by (iii), implies
that $x(\ker\,f)\# y(\ker\,f)$. Thus $\varphi$ is an se-mapping.

Let $x(\ker\,f)\# y(\ker\,f)$; that is, by (iii),
$(x,y)\in\mathrm{coker}\,f$. So we have $f(x)\# f(y)$, which, by the
definition of $\varphi$ means
$\varphi(x(\ker\,f))\#\varphi(y(\ker\,f))$. Thus $\varphi$ is
a-injective.

By the definition of composition of functions,
Theorem~\ref{ch-msa-t-coe}, and the definition of $\varphi
$,\textbf{\ }for each $x\in S$ we have
\[
(\varphi\pi)(x)=\varphi(\pi(x))=\varphi(x(\ker\,f))=f(x).
\]

(v)\ Taking into account (iv), we have to prove only that $\varphi$
is a surjection. Let $y\in T$. Then, as $f$ is onto, there exists
$x\in S$ such that $y=f(x)$. On the other hand $\pi(x)=x(\ker\,f)$.
By
(iv), we now have%
\[
y=f(x)=(\varphi\,\pi)(x)=\varphi(\pi(x))=\varphi(x(\ker\,f)).
\]
Thus $\varphi$ is a surjection.
\end{proof}

\begin{proof} \  -  \textbf{Theorem~\ref{ch-msa-BHM-1}}

(i). \ \ Let $x,y\in S$ and assume that $(x,y)\in
\varepsilon\cap\kappa$. Then $(x,y)\in \varepsilon$ and $(y,y)\in
\varepsilon$, which, by Lemma~\ref{ch-msa-l2b}, i.e. (Ap6') and and
$(x,y)\in\kappa$, gives $(y,y)\in \kappa$, which is impossible.
Thus, $\varepsilon\cap\kappa = \emptyset$.

Conversely, let $(x,a), (y,b)\in \varepsilon$ and $(x,y)\in\kappa$.
Then, by the co-transitivity of $\kappa$ and by the assumption, we
have
\begin{eqnarray*}
(x,y)\in \kappa &\ \Rightarrow\ & (x,a)\in \kappa \vee (a,y)\in
\kappa\\
 &\ \Rightarrow\ & (x,a)\in \kappa \vee (a,b)\in \kappa \vee (b,y)\in\kappa\\
 &\ \Rightarrow\ &  (a,b)\in \kappa.
\end{eqnarray*}

(ii). \ Let $\pi(x)\#\pi(y)$, that is $x\varepsilon \#
y\varepsilon$, which, by (i), means that $(x,y)\in \kappa$. Then, by
the strong irreflexivity of $\kappa$, we have $x\# y$. So $\pi$ is
an se-mapping.

Let $a\varepsilon\in S/\varepsilon$ and $x\in a\varepsilon$. Then
$(a,x)\in\varepsilon$, i.e. $a\varepsilon = x\varepsilon$, which
implies that $a\varepsilon = x\varepsilon=\pi(x)$. Thus $\pi$ is an
se-surjection.
\end{proof}

\begin{proof} \  -  \textbf{Theorem~\ref{ch-msa-basicfactor}}

(i). \  It follows from Theorem~\ref{ch-msa-BHM-1}(i).

(ii). \ It follows from Theorem~\ref{ch-msa-BHM-1}(ii).

(iii). \ This was shown in Theorem~\ref{ch-msa-cmrth2}(iii).

(iv). \ Let $\varphi$ be an se-mapping. Let $(x,y)\in
\mathrm{\mathrm{\mathrm{coker}}}\, f$ for some $x,y\in S$. Then, by
the definition of $\mathrm{coker}\,f$ and $\varphi$, the assumption
and (i), we have
\begin{eqnarray*}
f(x)\# f(y)  \ &\Leftrightarrow & \ \varphi(x(\mathrm{ker}\,f))\#
\varphi(y(\mathrm{ker}\,f))\\
&\Rightarrow  & \  x(\mathrm{ker}\,f) \# y(\mathrm{ker}\,f)\\
&\Leftrightarrow & \ (x, y)\in \kappa .
\end{eqnarray*}

Conversely, let $\mathrm{coker}\,f\subseteq\kappa$. By  assumption,
(i), and the definitions of $\varphi$ and $\mathrm{coker}\,f$, we
have
\begin{eqnarray*}
\varphi(x(\mathrm{ker}\,f))\# \varphi(y(\mathrm{ker}\,f)) \
&\Leftrightarrow & \ f(x)\# f(y)\\
&\Leftrightarrow & \ (x,y)\in \mathrm{coker}\,f\\
&\Rightarrow  & \ (x,y)\in \kappa\\
&\Leftrightarrow & \ x(\mathrm{ker}\,f)\# y(\mathrm{ker}\,f).
\end{eqnarray*}

(v). \ Let $\varphi$ be a-injective, and let $(x,y)\in \kappa$.
Then, by (i), we have
\begin{eqnarray*}
x(\mathrm{ker}\,f) \# y(\mathrm{ker}\,f) \ &\Rightarrow  & \
\varphi(x(\mathrm{ker}\,f))\# \varphi(y(\mathrm{ker}\,f))\\
  & \Leftrightarrow & \ f(x)\# f(y) \\
&\Leftrightarrow & \ (x,y)\in \mathrm{coker}\,f.
\end{eqnarray*}

Conversely, let $\kappa\subseteq \mathrm{coker}\,f$. Then
\begin{eqnarray*}
x(\mathrm{ker}\,f) \# y(\mathrm{ker}\,f) \  & \Leftrightarrow & \
(x,y)\in \kappa\\
&\Rightarrow  & \ (x,y)\in \mathrm{coker}\,f\\
&\Leftrightarrow & \ f(x)\# f(y) \\
&\Leftrightarrow & \ \varphi(x(\mathrm{ker}\,f))\#
\varphi(y(\mathrm{ker}\,f)).
\end{eqnarray*}

(vi). \ If $\varphi$ is an se-mapping then, by (iv), we have that
$\mathrm{coker}\,f\subseteq\kappa$. So, the strong irreflexivity of
$\kappa$ implies the  strong irreflexivity of $\mathrm{coker}\,f$,
i.e. $f$ is an se-mapping.
\end{proof}

\begin{proof} \  - \textbf{Lemma~\ref{ch-msa-cbl}}

(i).\ Follows by Lemma~\ref{ch-msa-sepl2}  and
Lemma~\ref{ch-msa-cun}.

(ii).\ This is a consequence of Theorem~\ref{ch-msa-t-coe}(ii).

(iii).\ Let $x\kappa^c_\tau, y\kappa^c_\tau\in S/\kappa^c_{\tau}$
such that $(x\kappa^c_\tau, y\kappa^c_\tau)\in\Upsilon_\tau$, i.e.
$(x,y)\in\tau$, which, by (ii), gives $x\kappa^c_\tau \#\,
y\kappa^c_\tau$. Thus, $\Upsilon_\tau$ is strongly irreflexive.

Let $x\kappa^c_\tau, y\kappa^c_\tau\in S/\kappa^c_{\tau}$ such that
$(x\kappa^c_\tau, y\kappa^c_\tau)\in\Upsilon_\tau$  and let
$z\kappa^c_\tau\in S/\kappa^c_{\tau}$. Then $(x,y)\in \tau$. By the
co-transitivity of $\tau$ we have $(x,z)\in \tau$ or $(z,y)\in
\tau$, which, in turn, gives $(x\kappa^c_\tau,
z\kappa^c_\tau)\in\Upsilon_\tau$ or $(z\kappa^c_\tau,
y\kappa^c_\tau)\in\Upsilon_\tau$. So, $\Upsilon_\tau$ is
co-transitive.

Let $x\kappa^c_\tau, y\kappa^c_\tau\in S/\kappa^c_{\tau}$ such that
$x\kappa^c_\tau \#\, y\kappa^c_\tau$. Then, by (ii), $(x,y)\in
\kappa_\tau= \tau\cup \tau^{-1}$. Thus $(x,y)\in \tau$ or $(y,x)\in
\tau$, which implies $(x\kappa^c_\tau,
y\kappa^c_\tau)\in\Upsilon_\tau$ or $(y\kappa^c_\tau,
x\kappa^c_\tau)\in\Upsilon_\tau$. Therefore, $\Upsilon_\tau$ is
co-antisymmetric. We have proved $\Upsilon_\tau$ is a co-order.

(iv). \ By Theorem~\ref{ch-msa-t-coe}(iii), $\pi_\tau$ is an
se-surjection. By the definition of $\Upsilon_\tau$, $\pi_\tau$ is
isotone and reverse isotone.
\end{proof}

\begin{proof} \  -  \textbf{Theorem~\ref{ch-msa-cp11}}

(i). \ Let $x,y\in S$ such that $(x,y)\in\mu_f$, that is
$(f(x),f(y))\in\sigma$. Now we have $f(x)\#\,f(y)$, and, finally, as
$f$ is an se-mapping, $x\#\,y$ follows. Thus, $\mu_f$ is strongly
irreflexive.

Let $x,y\in S$ such that $(x,y)\in\mu_f$, and let $z\in S$. Then, by
the definition of $\mu_f$, $(f(x),f(y))\in\sigma$. By
co-transitivity of $\sigma$ we have $(f(x),f(z))\in\sigma$ or
$(f(z),f(y))\in\sigma$, that is $(x,z)\in\mu_f$ or $(z,y)\in\mu_f$.
So, $\mu_f$ is co-transitive.

(ii).\ By Lemma~\ref{ch-msa-cbl}(i) $\kappa_f=\mu_f\cup\mu_f^{-1}$
is a co-equivalence on on $S$. Let $(x,y)\in\kappa_f$, i.e.
$(x,y)\in \mu_f\cup\mu_f^{-1}$, that is  $(x,y)\in \mu_f$ or
$(y,x)\in \mu_f$. So we have $(f(x),f(y))\in\sigma$ or
$(f(y),f(x))\in\sigma$. Strong irreflexivity of $\sigma$ implies
$f(x)\#\,f(y)$, i.e. $(x,y)\in \mathrm{coker}\,f$. Thus,
$\kappa_f\subseteq \mathrm{coker}\,f$.

(iii).\ Let $x,y\in S$ such that $(x,y)\in \mathrm{coker}\,f$, i.e.
$f(x)\#\,f(y)$. Now, by the co-antisymmetry of $\sigma$, we have
$(f(x),f(y))\in\sigma$ or $(f(y),f(x))\in\sigma$, that is, by the
definition of $\mu_f$, $(x,y)\in \mu_f$ or  $(y,x)\in \mu_f$, i.e.
$(x,y)\in \mu_f\cup\mu_f^{-1} = \kappa_f$.

(iv).\ Follows immediately by Lemma~\ref{ch-msa-cbl}(ii),(iii).

(v).\ Let us, first, prove that $\psi$ is well-defined. Let
$x\kappa_f = y\kappa_f$. Then
\begin{eqnarray*}
x\kappa^c_f = y\kappa^c_f &\Leftrightarrow & (x,y)\in \kappa^c_f\\
&\Leftrightarrow & \neg ((x,y)\in \kappa_f)\\
&\Leftrightarrow & \neg ((x,y)\in \mu_f \vee (y,x)\in\mu_f)\\
&\Rightarrow & \neg ((f(x),f(y))\in \sigma \vee (f(y),f(x))\in \sigma)\\
&\Leftrightarrow & \neg ((f(x),f(y))\in \sigma \cup \sigma^{-1}=
\kappa_\sigma)\\
&\Leftrightarrow & (f(x),f(y))\in \kappa^c_\sigma\\
&\Leftrightarrow & (f(x))\kappa^c_\sigma = (f(y))\kappa^c_\sigma\\
&\Leftrightarrow & \psi( x\kappa^c_f) = \psi (y\kappa^c_f).
\end{eqnarray*}
Let $\psi( x\kappa^c_f) \#\, \psi (y\kappa^c_f)$. Then
\begin{eqnarray*}
\psi( x\kappa^c_f) \#\, \psi (y\kappa^c_f)&\Leftrightarrow &
(f(x))\kappa^c_\sigma \#\, (f(y))\kappa^c_\sigma\\
&\Leftrightarrow & (f(x),f(y))\in \kappa^c_\sigma \\
&\Leftrightarrow & (f(x),f(y))\in \sigma \vee (f(y),f(x))\in
\sigma\\
&\Rightarrow & (x,y)\in \mu_f \vee (y,x)\in\mu_f \\
&\Leftrightarrow & (x,y)\in\kappa_f\\
&\Leftrightarrow & x\kappa^c_f \#\, y\kappa^c_f,
\end{eqnarray*}
so $\psi$ is an se-mapping.

Let $x\in S$. Then
$$(\psi\pi_f )(x)= \psi(\pi_f (x))= \psi(x\kappa^c_f)=
(f(x))\kappa^c_\sigma = \pi_\sigma (f(x))= (\pi_\sigma f) (x).$$
\end{proof}

\begin{proof} \  -  \textbf{Theorem~\ref{ch-msa-cqo12}}

Let us, first prove that $\varphi$ is well-defined. Let
$x\kappa^c_\tau, y\kappa^c_\tau\in S/\kappa^c_{\tau}$ be such that
$x\kappa^c_\tau = y\kappa^c_\tau$, i.e. $(x,y)\bowtie \kappa_\tau$.
Assume $\varphi(x\kappa^c_\tau)\#\, \varphi(y\kappa^c_\tau)$, i.e.
by the definition of $\varphi$, $f(x)\,\#\, f(y)$, that is
$(x,y)\in\mathrm{coker}\,f$. Now, by the
Theorem~\ref{ch-msa-cp11}(iii), we have $(x,y)\in \kappa_f$ which is
a contradiction. Thus, $\neg(f(x)\,\#\, f(y))$, which,as apartness
on $T$ is tight one, implies $f(x)= f(y)$.

Let $x\kappa^c_\tau, y\kappa^c_\tau\in S/\kappa^c_{\tau}$ be such
that $\varphi(x\kappa^c_\tau)\,\#\, \varphi(y\kappa^c_\tau)$, i.e.
$f(x)\,\#\, f(y)$. Therefore, $(x,y)\in\mathrm{coker}\,f =
\kappa_f$. From the assumption $\mu_f \subseteq \tau$, it follows
that $\kappa_f\subseteq \kappa_\tau$, which, in turn, implies
$(x,y)\in \kappa_\tau$. Now, by Lemma~\ref{ch-msa-cbl}(ii),
$x\kappa^c_\tau\,\#\, y\kappa^c_\tau$, and $\varphi$ is an
se-mapping.

Finally, let $x\in S$. Then
$$\varphi\pi_\tau (x)=\varphi(\pi_\tau (x))=
\varphi(x\kappa^c_\tau)= f(x).$$
\end{proof}

\begin{proof} \  - \textbf{Theorem~\ref{ch-msa-evid}}

As  in the classical case, composition of functions is associative.
By Theorem~\ref{ch-msa-thm1.2}, $(\mathcal{T}_X^{se}, =, \#)$ is a
set with apartness.

Let $f, g\in \mathcal{T}_X^{se}$ and suppose that $(fg)(x)\# (f
g)(y)$ for some $x,y\in X$. Then, by the definition of the
composition, $f(g(x))\#\, f(g(y))$, and, as $f$ is an se-mapping, we
have $g(x)\#g(y)$. Finally, as $g$ is an se-mapping as well, we have
$x\#y$. Thus, $f g$ is an se-mapping and $f g\in
\mathcal{T}_X^{se}$.

Let $f,g,h,w\in \mathcal{T}_X^{se}$ and $f h\# g w$. Then, by the
definition of apartness in $S$, there is an element $x\in X$ such
that $(f h)(x)\#(g w)(x)$, i.e. $f(h(x))\#g(w(x))$. Now we have
$$f(h(x))\# f(w(x)) \vee f(w(x))\#\, g(w(x)),$$
which, further, implies $h(x)\#w(x)$ (because $f$ is an se-mapping)
or $f\#g$ (by the definition of the apartness relation on
$\mathcal{T}_X^{se}$). Thus $f\#g \vee h\#w,$ that is, composition
is an se-operation and $(\mathcal{T}_X^{se}, =, \# )$ is a semigroup
with apartness.
\end{proof}

\begin{proof} \  - \textbf{Theorem~\ref{ch-msa-Cayley}}

Let $(S,= , \# ;\,  \cdot)$  be a semigroup with apartness. The
semigroup $S$ embeds into the monoid with apartness
$S^1=(S\cup\{1\},=_1, \# _1, \cdot)$ with equality $=_1$ which
consists of all pairs in $=$ and the pair $(1,1)$, and with
apartness $\# _1$ which consists of all pairs in $\# $ and the pairs
$(a,1),(1,a)$ for each $a\in S$.

Let $f_a$ be a left translation of $S^1$, i.e. $f_a(x) = ax$, for
all $x\in S^1$. Then $f_a$ is an se-function.  Indeed, $f_a(x)\#_1
f_a(y)$  is equivalent to $ax \#_1 ay$. The strong extensionality of
multiplication implies $x \#_1  y$.

Denote by $\mathcal{T}_{S^1}^{se}$ the set of all se-functions from
$S^1$ to $S^1$.  In a pretty much similar way as in \textbf{CLASS}
define a mapping $\varphi: S^1\rightarrow \mathcal{T}_{S^1}^{se}$ to
be $ \varphi(a) = f_a$, for each $a\in S^1$. It is routine to verify
that
$$
\varphi(ab) =  f_{ab} =  f_a\,f_b = \varphi(a)\varphi(b),
$$
as well as
$$
\varphi(a) \# _T \varphi(b) \ \Rightarrow \ a\#_1 b.
$$
Thus, $\varphi$ is an se-homomorphism. Also, $\varphi(a)=_T
\varphi(b)$ iff  $ax = _1 bx$ for all $x\in S^1$, and, for $x=1$, we
have $a=_1 b$. Therefore, $\varphi$ is an embedding.
\end{proof}

\begin{proof} \  -  \textbf{Lemma~\ref{ch-msa-lrc}}

Let $\tau$ be a co-compatible co-quasiorder w.r.f. multiplication on
$S$, and let $x,y,z\in S$. Then $(zx,zy)\in \tau$ implies
$(x,y)\in\tau$ or $(z,z)\in\tau$. The latter is impossible because
of strong irreflexivity of  $\tau$. Thus $(x,y)\in\tau$, i.e. $\tau$
is left co-compatible.

Conversely, let $\tau$ be a left and a right co-compatible
co-quasiorder on $S$. Let $x,y,z,t\in S$ be such that $(xz,yt)\in
\tau$. By the co-transitivity of $\tau$, it follows either
$(xz,yz)\in \tau$ or $(yz,yt)\in \tau$. Now, by the assumption, we
have $(x,y)\in \tau$ or $(z,t)\in\tau$, as required.
\end{proof}

\begin{proof} \  - \textbf{Theorem~\ref{ch-msa-sap2}}

\ By Theorem~\ref{ch-msa-t-coe}, $(S/\kappa^c,=,\neq)$ is a set with
apartness. The associativity of multiplication in $S/\kappa^c$
follows from that of multiplication on $S$.

Let $a\kappa^c\,x\kappa^c\# b\kappa^c\,y\kappa^c$. Then
$(ax)\kappa^c\# (by)\kappa^c$. By Theorem~\ref{ch-msa-t-coe}, we
have that $(ax,by)\in\kappa$. But $\kappa$ is a co-congruence, so
either $(a,b)\in\kappa$ or $(x,y)\in\kappa$. Thus, by the definition
of $\#$ in $S/\kappa^c$, either $a\kappa^c\# b\kappa^c$ or
$x\kappa^c\# y\kappa^c$. So $(S/\kappa^c ,=,\#,\cdot\,)$ is a
semigroup with apartness. Using that fact and the definition of
$\pi$, we have
\[
\pi(xy)=(xy)\kappa^c=x\kappa^c\,y\kappa^c=\pi(x)\pi(y).
\]
Hence $\pi$ is a homomorphism, and, by Theorem~\ref{ch-msa-t-coe} ,
$\pi$ is an se-surjection.
\end{proof}

\begin{proof} \  -  \textbf{Theorem~\ref{ch-msa-cmrth22}}

$\ $(i) \ Taking into account Theorem~\ref{ch-msa-cmrth2}, it is
enough to prove that $\mathrm{coker}\,f$ is co-compatible with
multiplication in $S$. Let $(ax,by)\in\mathrm{coker}\,f$, i.e.
$f(ax)\# f(by)$. Since $f$ is a homomorphism, we have $f(a)f(x)\#
f(b)f(y)$. The strong extensionality of multiplication implies that
either $f(a)\# f(b)$ or $f(x)\# f(y)$. Thus either $(a,b)\in
\mathrm{coker}\,f$ or $(x,y)\in\mathrm{coker}\,f$, and therefore
$\mathrm{coker}\,f$ is a co-congruence on $S$.

(ii). \ This is Theorem~\ref{ch-msa-cmrth2}(ii).

(iii) \ This follows by Theorem~\ref{ch-msa-cmrth2} and
Theorem~\ref{ch-msa-sap2}.

(iv) \ Using (iii) and the assumption that $f$ is a homomorphism, we
have
\begin{align*}
\varphi(x(ker\,f)\,y(ker\,f))  & =\varphi((xy)(ker\,f))\\
& =f(xy)\\
& =f(x)f(y)\\
& =\varphi(x(ker\,f))\,\varphi(y(ker\,f)).
\end{align*}\t
Now, by Theorem~\ref{ch-msa-cmrth2}, $\varphi$ is  an apartness
embedding.

(v) \ This follows by Theorem~\ref{ch-msa-cmrth2} and (iv).
\end{proof}

\begin{proof} \  -  \textbf{Theorem~\ref{ch-msa-cmrth3}}

(i). \ If $\kappa$ defines an apartness on $S/\mu$, then, by
Theorem~\ref{ch-msa-BHM-1}(i), $\mu\cap\kappa=\emptyset$.

Let $\mu$ be a congruence and $\kappa$ a co-congruence on a
semigroup with apartness $S$ such that $\mu\cap\kappa=\emptyset$.
Then, by Theorem~\ref{ch-msa-BHM-1}(i), $\kappa$ defines apartness
on $S/\mu$.

Let $a\mu\,x\mu \# \,b\mu \, y\mu$, then $(ax)\mu\#\,(by)\mu$ which
further, by the definition of apartness on $S/\mu$, ensures that
$(ax,by)\in\kappa$. But $\kappa$ is a co-congruence, so either
$(a,b)\in\kappa$ or $(x,y)\in\kappa$. Thus, by the definition of
apartness in $S/\mu$ again, either $a\mu \# b\mu$ or $x\mu \# y\mu$.
So $(S/\mu ,=,\#,\cdot\,)$ is a semigroup with apartness.

(ii).\ Straightforward.
\end{proof}

\begin{proof} \  - \textbf{Lemma~\ref{ch-msa-cb2}}

(i).\ By Lemma~\ref{ch-msa-cbl}(i), $\kappa_{\tau}$ is a
co-equivalence. Let $(ax,by)\in \kappa_\tau$. Then, by the
definition of $\kappa_\tau$, we have $(ax,by)\in \tau$ or
$(by,ax)\in \tau$. By the co-compability of $\tau$ we have
$(a,b)\in\tau$ or $(x,y)\in\tau$ or $(b,a)\in\tau$ or
$(x,y)\in\tau$, i.e. $(a,b)\in \kappa_\tau$ or $(x,y)\in
\kappa_\tau$. So, $\kappa_\tau$ is co-compatible, and, therefore,
co-congruence.

(ii).\ Follows by (i) and Theorem~\ref{ch-msa-sap2}.

(iii).\ By Lemma~\ref{ch-msa-cbl}(ii), it follows that
$\Upsilon_\tau$ is a co-order. Let us prove that it is
co-compatible. Let $x\kappa^c_\tau, y\kappa^c_\tau, a\kappa^c_\tau,
b\kappa^c_\tau\in S/\kappa^c_\tau$ be such that
$(x\kappa^c_\tau\,a\kappa^c_\tau , y\kappa^c_\tau\,
b\kappa^c_\tau)\in \Upsilon_\tau$. Then $((xa)\kappa^c_\tau ,
(yb)\kappa^c_\tau)\in \Upsilon_\tau$, i.e. $(xa,yb)\in\tau$. By the
co-compatibility of $\tau$, we now have $(x,y)\in\tau$ or
$(a,b)\in\tau$. By the definition of $\Upsilon_\tau$, it follows
that $(x\kappa^c_\tau , y\kappa^c_\tau)\in \Upsilon_\tau$ or
$(a\kappa^c_\tau , b\kappa^c_\tau)\in \Upsilon_\tau$; that is
$\Upsilon_\tau$ is co-compatible.

(iv).\ Follows by Theorem~\ref{ch-msa-sap2} and
Lemma~\ref{ch-msa-cbl}(iii).
\end{proof}

\end{document}